\documentclass[10pt]{article}       
\title{The tree method for multidimensional $q$-Hahn and $q$-Racah polynomials }

\author{Fabio Scarabotti}

\usepackage{latexsym,amsbsy,amsmath,amscd,amssymb}

\usepackage{latexsym,amsbsy,amsmath,amscd,amssymb,epic,eepicemu}  

\usepackage{amsthm}       
\usepackage{amssymb,amsmath,latexsym}
 \newtheorem{definition}{Definition} [section]       
 \newtheorem{remark}[definition]{Remark}       
 \newtheorem{example}[definition]{Example}

 \newtheorem{proposition}[definition]{Proposition}       
 \newtheorem{theorem}[definition]{Theorem}       
 \newtheorem{corollary}[definition]{Corollary}       
  \newtheorem{lemma}[definition]{Lemma}

\begin{document}

\maketitle

\begin{abstract}  

We develop a tree method for multidimensional $q$-Hahn polynomials. We define them as eigenfunctions of a multidimensional $q$-difference operator and we use the factorization of this operator as a key tool. Then we define multidimensional $q$-Racah polynomials as the connection coefficients between different bases of $q$-Hahn polynomials. We show that our multidimensional $q$-Racah polynomials may be expressed as product of ordinary one-dimensional $q$-Racah polynomial by means of a suitable sequence of transplantations of edges of the trees. Our paper is inspired to the classical tree methods in the theory of Clebsch-Gordan coefficients and  of hyperspherical coordinates. It is based on previous work of Dunkl, who considered two-dimensional $q$-Hahn polynomials. It is also related to a recent paper of Gasper and Rahman: we show that their multidimensional $q$-Racah polynomials correspond to a particular case of our construction.\\



\footnote{FABIO SCARABOTTI, Dipartimento MeMoMat, Universit\`a degli Studi di Roma ``La Sapienza'', via A. Scarpa 8, 00161 Roma (Italy)
{\it e-mail:} {\tt scarabot@dmmm.uniroma1.it}\\
}
\end{abstract}

%%%%%%%%%%%%%%%%%%%%%%%%%%%%%%%%%%%%%%%%%%%%%%%%%%%%%%%%%%%%%%%%%%%%%%%%%%%%%%%%%%
%%%%%%%%%%%%%%%%%%%%%%%%%%%%%%%%%%%%%%%%%%%%%%%%%%%%%%%%%%%%%%%%%%%%%%%%%%%%%%%%%%%%%
\section{Introduction}
%%%%%%%%%%%%%%%%%%%%%%%%%%%%%%%%%%%%%%%%%%%%%%%%%%%%%%%%%%%%%%%%%%%%%%%%%%%%%%%%%%
%%%%%%%%%%%%%%%%%%%%%%%%%%%%%%%%%%%%%%%%%%%%%%%%%%%%%%%%%%%%%%%%%%%%%%%%%%%%%%%%%%%%%

The tree method is a powerful procedure to construct systems of multidimensional orthogonal polynomials that are expressible in terms orthogonal polynomials in one variable. Its origin is in the theory of Clebsch-Gordan coefficients \cite{BL,JLV,NSU,VdJ} and in the theory of hyperspherical coordinates \cite{K-V, NSU}. In a previous paper \cite{ScarabottimultHahn} we used this method to construct multidimensional Hahn polynomials as intertwining functions on the symmetric group, extending some results of Dunkl \cite{Du5}. The aim of the present paper is to develop a tree method for multidimensional $q$-Hahn polynomials and to show that the connection coefficients between different bases of multidimensional $q$-Hahn polynomials can be naturally interpreted as multidimensional $q$-Racah polynomials. In this way we extend Dunkl's results in the paper \cite{Du6}.\\

A family of multidimensional $q$-Racah polynomials has been recently introduced also by Gasper and Rahman \cite{GaRa}, generalizing previous results of Tratnik \cite{Tratnik} in the $q=1$ case; they also discuss $q$-Hahn and other polynomials as limiting cases. We show that their polynomials may be obtained as particular cases of our construction. More precisely, their multidimensional $q$-Hahn polynomials correspond to the polynomials that we associate to a tree in which every right subtree has only one branch. Their multidimensional $q$-Racah polynomials correspond to the connection coefficients between the multidimensional $q$-Hahn polynomials associated respectively to the binary tree in 
which every left subtree has only one branch and to the binary tree in which every right
subtree has only one branch.\\

In contrast to \cite{ScarabottimultHahn}, in the present paper we do not use group theoretical methods nor we give group theoretical interpretations of our results. We use a purely analytic approach, mainly based on the manipulation of finite $q$-difference operators. Indeed, in the $q$-setting a group theoretical approach should be based on the representation theory of the finite general linear group, which is not developed nor manageable as the representation theory of the symmetric group (that we used in \cite{ScarabottimultHahn}); see \cite{CST3,DuqHahn,MaPa} for the one dimensional case. Another possible group theoretical approach to multidimensional $q$-Hahn polynomials is through the representation theory of quantum groups: see \cite{Rosen}.\\

This is the plan of the paper. 
In Section \ref{secmultqHahnop} we introduce a multidimensional $q$-Hahn operator and derive its spectral analysis by means of a suitable factorization method. 
In Section \ref{secqHahn} we derive the basic properties of the one-dimensional $q$-Hahn polynomials using the factorization method. 
In Section \ref{building} we use the factorization method and some properties of the $q$-Hahn polynomials to establish a basic recursive procedure. In Section \ref{treemethod} we use such recursive procedure to define the multidimensional $q$-Hahn polynomials and to derive their basic properties. We also give three examples, one of which coincides with the $q$-Hahn polynomials in \cite{GaRa}. In Section \ref{bidqHahn} we reproduce Dunkl's proof of the fact that the $q$-Racah polynomials are connection coefficients between two bases of bidimensional $q$-Hahn polynomials. In Section \ref{sectiontransplantation} we show that the result in the previous section may be used to perform the basic operation of transplantation of an edge, that is to express the $q$-Hahn polynomials associated to the tree $\mathcal{T}$ in terms of the polynomials associated to the tree $\mathcal{S}$ if $\mathcal{T}$ is obtained from $\mathcal{S}$ by means of the transplantation of an edge. Actually, we have arranged the normalization constants in order to get the formulas for transplantations from right to left as simple as possible at the cost that the coefficients for a transplantation from left to right are more complicated. In Section \ref{sectionmultqRacah} we define the multidimensional $q$-Racah polynomials as the connection coefficients between bases of multidimensional $q$-Hahn polynomials associated to different trees. We also show that the result in the previous section may be used to obtain explicit formulas for the multidimensional polynomials in terms of one dimensional $q$-Racah polynomials when the trees may be connected by a sequence of right transplantations of edges. In Section \ref{sectionGaRa} we analyze the multidimensional $q$-Racah polynomials in \cite{GaRa}. \\

Our paper is quite self-contained since our construction naturally gives many properties of $q$-Hahn and $q$-Racah polynomials even in the one-dimensional case. Just in two cases, to get the standard $\;_3\phi_2$ (resp. $\;_4\phi_3$) expression for the $q$-Hahn (resp. $q$-Racah) polynomials, we will use nontrivial identities from the theory of basic hypergeometric series. We refer to \cite{GaRabook} for those identities and also for unexplained terminology on $q$-shifted factorials (but we will use nonstandard normalizations for one-dimensional $q$-Hahn and $q$-Racah polynomials).

%%%%%%%%%%%%%%%%%%%%%%%%%%%%%%%%%%%%%%%%%%%%%%%%%%%%%%%%%%%%%%%%%%%%%%%%%%%%%%%%%%
%%%%%%%%%%%%%%%%%%%%%%%%%%%%%%%%%%%%%%%%%%%%%%%%%%%%%%%%%%%%%%%%%%%%%%%%%%%%%%%%%%%%%
\section{The multidimensional $q$-Hahn operator and its factorization}\label{secmultqHahnop}
%%%%%%%%%%%%%%%%%%%%%%%%%%%%%%%%%%%%%%%%%%%%%%%%%%%%%%%%%%%%%%%%%%%%%%%%%%%%%%%%%%%%%

We introduce here several notation that will be used throughout the paper.
We use boldface letters to denote ordered sequences of nonnegative integers (the {\em variables}): ${\bf x}=(x_1,x_2,\dotsc,x_h)$ means that each $x_i$ is a nonnegative integer. If $N$ is a positive integer, an $h$-{\em parts} composition of $N$ is an ordered sequence ${\bf x}=(x_1,x_2,\dotsc,x_h)$ of nonnegative integers such that $x_1+x_2+\dotsb+x_h=N$. We denote by
$[h;N]$ the set of all $h$-parts compositions of $N$ and by $V_{h,N}$ the vector space of all complex valued functions defined on $[h;N]$.
We will also use two specific notation: we set $\epsilon^\pm_j(x_1,\dotsc,x_j,\dotsc,x_h)=(x_1,\dotsc,x_j\pm 1,\dotsc,x_h)$ and, for a couple of positive integers $i,j$, we set $[i<j]=1$ if $i<j$, $[i<j]=0$ if $i\geq j$.
If $f\in V_{h,N}$, ${\bf x}\in [h;N]$ and $x_j=N$ (or $x_i=0$) we set $f(\epsilon_j^+({\bf x}))=0$ (respectively $f(\epsilon_i^-({\bf x}))=0$).\\

Now fix a real number $0<q<1$, a positive integer $h$ and real numbers $\alpha_1,\alpha_2,\dotsc,\alpha_h$ (the {\em parameters}) satisfying the conditions $0<\alpha_i<q^{-1}$, $i=1,2,\dotsc,h$, or the conditions $\alpha_i>q^{-N}$, $i=1,2,\dotsc,h$. We will use the following notation for the variables and the parameters: we set

\[
X_k=\sum_{j=1}^kx_j \qquad\quad\text{and}\qquad\quad A_k=\prod_{j=1}^ka_j, \qquad\quad\text{for } k=1,2,\dotsc,h.
\]

\noindent
We introduce a scalar product on $V_{h,N}$ by setting

\begin{equation}\label{scalarhmult}
\langle f_1,f_2 \rangle_{V_{h,N}}=q^{N(N+1)/2}\sum_{{\bf x}\in [h;N]}\prod_{i=1}^h\left[\frac{(q\alpha_i;q)_{x_i}}{(q;q)_{x_i}}\left(\alpha_i q\right)^{N-X_i}\right]f_1({\bf x})\overline{f_2({\bf x})}.
\end{equation}

\noindent
It is easy to see that if we set $h=s+1$ and $\alpha_i=a_i$, $i=1,2,\dotsc,s+1$, the weight in \eqref{scalarhmult} coincides with (3.17) in \cite{GaRa} multiplied by the factor 

\[
\frac{(a_{s+1}q;q)_N}{(q;q)_N}(A_sq^s)^Nq^{N(N+1)/2}.
\]

The {\em multidimensional} $q$-{\em Hahn operator} $\mathcal{D}_N:V_{h,N}\rightarrow V_{h,N}$ is defined by setting, for $f\in V_{h,N}$,

\begin{multline}\label{multqHahnoperator}
\mathcal{D}_Nf({\bf x})=\sum_{\substack{i,j=1\\i\neq j}}^hA_{j-1}q^{j-1+X_{j-1}+X_{i-1}-N-[i<j]}(\alpha_jq^{x_j+1}-1)(1-q^{x_i})
f(\epsilon_i^-\epsilon_j^+({\bf x}))\\
+\biggl\{\sum_{j=1}^hA_{j-1}q^{j-1+2X_{j-1}-N}(\alpha_jq^{x_j}-1)(1-q^{x_j})
+(A_h q^{h+N-1}-1)(1-q^{-N})\biggr\}f({\bf x}),
\end{multline}

\noindent 
for all ${\bf x}\in[h;N]$. We associate to $\mathcal{D}_N$ a {\em raising operator} $\mathcal{R}_N:V_{h,N}\rightarrow V_{h,N+1}$ and a {\em lowering operator} $\mathcal{L}_N:V_{h,N}\rightarrow V_{h,N-1}$ defined by setting, for $f\in V_{h,N}$, respectively

\[
\begin{aligned}
&\mathcal{R}_Nf({\bf x})=\sum_{i=1}^hq^{X_{i-1}-N-1}(1-q^{x_i})f(\epsilon_i^-({\bf x})),\qquad\qquad\qquad\qquad\qquad\quad&\text{for all }\quad{\bf x}\in[h;N+1],\\
&\mathcal{L}_Nf({\bf x})=\sum_{j=1}^hA_{j-1}q^{j-1+X_{j-1}}(\alpha_jq^{x_j+1}-1)f(\epsilon_j^+({\bf x})),\qquad&\text{for all }\quad{\bf x}\in[h;N-1].
\end{aligned}
\]

\noindent
We also denote by $I_N$ the identity operator on $V_{h,N}$. Now we prove a series of identities for those operators. At the end we will get the spectral analysis of $\mathcal{D}_N$.

\begin{lemma}\label{multcomrel}
The raising and lowering operators satisfy the identities:

\begin{equation}\label{RN1LN}
\mathcal{R}_{N-1}\mathcal{L}_N=\mathcal{D}_N-(1-q^{-N})(A_h q^{N+h-1}-1)I_N,\qquad\qquad N\geq 1,
\end{equation}

\begin{equation}\label{LN1RN}
\mathcal{L}_{N+1}\mathcal{R}_N=\mathcal{D}_N-(1-q^{-N-1})(A_h q^{N+h}-1)I_N, \qquad\qquad N\geq 0,
\end{equation}

\noindent
and the commutation relation:

\begin{equation}\label{commrelmult}
\mathcal{L}_{N+1}\mathcal{R}_N-\mathcal{R}_{N-1}\mathcal{L}_N=q^{-N-1}(1-q)(A_h q^{2N+h}-1)I_N.
\end{equation}

\end{lemma}

\begin{proof}
If $f\in V_{h,N}$ and ${\bf x}=(x_1,\dotsc,x_h)\in [h;N]$ we have:

\[
\begin{split}
\mathcal{R}_{N-1}\mathcal{L}_Nf({\bf x})=&\sum_{i=1}^hq^{X_{i-1}-N}(1-q^{x_i})\mathcal{L}_Nf(\epsilon_i^-({\bf x}))\\
=&\sum_{\substack{i,j=1\\j\neq i}}^hA_{j-1}q^{j-1+X_{j-1}-[i<j]}(\alpha_jq^{x_j+1}-1)
 q^{X_{i-1}-N}(1-q^{x_i})f(\epsilon_i^-\epsilon_j^+({\bf x}))\\
&+\biggl\{\sum_{j=1}^hA_{j-1}q^{j-1+2X_{j-1}-N}(\alpha_jq^{x_j}-1)(1-q^{x_j})\biggr\}f({\bf x})\\
=&\mathcal{D}_Nf({\bf x})-(A_h q^{N+h-1}-1)(1-q^{-N})f({\bf x}),
\end{split}
\]

\noindent
and therefore \eqref{RN1LN} is proved (the terms $-[i<j]$ appears after the second equality because if $i<j$ then $X_{j-1}$ must be replaced by $X_{j-1}-1$).
The proof of the \eqref{LN1RN} is similar but slightly more complicated: now we have

\[
\begin{split}
\mathcal{L}_{N+1}\mathcal{R}_N-\mathcal{D}_N=&\biggl\{\sum_{j=1}^hA_{j-1}q^{j-2+2X_{j-1}-N}
\bigl[(\alpha_jq^{x_j+1}-1)(1-q^{x_j+1})-q(\alpha_jq^{x_j}-1)(1-q^{x_j})\bigr]\biggr\}I_N\\
&-(A_h q^{N+h-1}-1)(1-q^{-N})I_N\\
=&-\sum_{j=1}^hA_jq^{j+2X_j-N}
-\sum_{j=1}^hA_{j-1}q^{j-2+2X_{j-1}-N}+\sum_{j=1}^hA_jq^{j-1+2X_j-N}\\
&+\sum_{j=1}^hA_{j-1}q^{j-1+2X_{j-1}-N}-(A_h q^{N+h-1}-1)(1-q^{-N})I_N\\
=&\left[-A_hq^{N+h}-q^{-N-1}+A_hq^{N+h-1}+q^{-N}-(A_h q^{N+h-1}-1)(1-q^{-N})\right]I_N\\
=&-(1-q^{-N-1})(A_h q^{N+h}-1)I_N.
\end{split}
\]

\noindent
In third equality most of the terms in the four summations cancel each other giving $-A_hq^{N+h}-q^{-N-1}+A_hq^{N+h-1}+q^{-N}$ as the final result. Finally, \eqref{commrelmult} follows immediately from \eqref{RN1LN} and \eqref{LN1RN}.

\end{proof}

\begin{lemma}

The operator $-\mathcal{L}_N$ is the adjoint of $\mathcal{R}_{N-1}$, that is

\begin{equation}\label{invscalarprod}
\langle\mathcal{L}_Nf_1,f_2\rangle_{V_{h,N-1}}=-\langle f_1,\mathcal{R}_{N-1}f_2\rangle_{V_{h,N}},
\end{equation} 

\noindent
for all $f_1\in V_{h,N}$ and $f_2\in V_{h,N-1}$. The operator $\mathcal{D}_N$ is self-adjoint. 

\end{lemma}

\begin{proof}

We have

\[
\begin{split}
\langle\mathcal{L}_Nf_1,f_2\rangle_{V_{h,N-1}}=&\sum_{{\bf x}\in [h;N-1]}\sum_{j=1}^hA_{j-1}q^{j-1+X_{j-1}}(\alpha_jq^{x_j+1}-1)\\
&\times q^{(N-1)N/2}\prod_{i=1}^h\left[\frac{(q\alpha_i;q)_{x_i}}{(q;q)_{x_i}}\left(\alpha_i q\right)^{N-1-X_i}\right]f_1(\epsilon_j^+({\bf x}))\overline{f_2({\bf x})}\\
({\bf y}=\epsilon^+_j({\bf x}),\quad y_j=x_j+1)\quad\qquad=&-\sum_{{\bf y}\in [h;N]}\sum_{j=1}^hq^{Y_{j-1}-N}(1-q^{y_j})\\
(\frac{(q\alpha_i;q)_{x_j}}{(q;q)_{x_j}}=\frac{(q\alpha_i;q)_{y_j}}{(q;q)_{y_j}}\frac{1-q^{y_j}}{1-\alpha_jq^{y_j}})\quad\qquad&\times q^{N(N+1)/2}\prod_{i=1}^h\left[\frac{(q\alpha_i;q)_{y_i}}{(q;q)_{y_i}}\left(\alpha_i q\right)^{N-Y_i}\right]f_1({\bf y})\overline{f_2(\epsilon^-_j{\bf y})}\\
=&-\langle f_1,\mathcal{R}_{N-1}f_2\rangle_{V_{h,N}}.
\end{split}
\]

\noindent
The selfadjointness of $\mathcal{D}_N$ may be obtained as a consequence of \eqref{invscalarprod} and \eqref{RN1LN}.

\end{proof}

\begin{lemma}
For $N>n\geq0$, we have:

\begin{equation}\label{LRRR}
\begin{split}
\mathcal{L}_N\mathcal{R}_{N-1}\mathcal{R}_{N-2}\dotsb \mathcal{R}_n=&\mathcal{R}_{N-2}\mathcal{R}_{N-3}\dotsb \mathcal{R}_{n-1}\mathcal{L}_n\\
&+q^{-N}(1-q^{N-n})(A_h q^{N+n+h-1}-1)\mathcal{R}_{N-2}\mathcal{R}_{N-3}\dotsb \mathcal{R}_n.
\end{split}
\end{equation}

\end{lemma}

\begin{proof}
For $N=n+1$ it coincides with \eqref{commrelmult}. The general case can be proved by induction on $N-n$:

\[
\begin{split}
\mathcal{L}_N\mathcal{R}_{N-1}\mathcal{R}_{N-2}\dotsb \mathcal{R}_n=&\mathcal{R}_{N-2}\mathcal{L}_{N-1}\mathcal{R}_{N-2}\mathcal{R}_{N-3}\dotsb \mathcal{R}_n
+q^{-N}(1-q)(A_h q^{2N+h-2}-1)\mathcal{R}_{N-2}\mathcal{R}_{N-3}\dotsb \mathcal{R}_n\\
=&\mathcal{R}_{N-2}\mathcal{R}_{N-3}\dotsb \mathcal{R}_{n-1}\mathcal{L}_n
+q^{-N}(1-q^{N-n})(A_h q^{N+n+h-1}-1)\mathcal{R}_{N-2}\mathcal{R}_{N-3}\dotsb \mathcal{R}_n,
\end{split}
\]

\noindent
where the first equality follows from \eqref{commrelmult} and the second equality from the inductive hypothesis (and an elementary algebraic calculation).
\end{proof}

\begin{lemma}
If $0\leq n\leq m\leq N$ and $\mathcal{L}_nf=0$ then

\begin{multline}\label{LLLRRR}
\mathcal{L}_{m+1}\mathcal{L}_{m+2}\dotsb\mathcal{L}_N\mathcal{R}_{N-1}\mathcal{R}_{N-2}\dotsb \mathcal{R}_nf\\
=(-1)^{N-m}q^{-(N-m)(N+m+1)/2}(q^{m-n+1};q)_{N-m}(A_h q^{n+m+h};q)_{N-m}\mathcal{R}_{m-1}\mathcal{R}_{m-2}\dotsb \mathcal{R}_nf.
\end{multline}

\end{lemma}
\begin{proof}
For $m=N-1$ it is an immediate consequence of \eqref{LRRR}; the general case follows by induction on $N-m$: if $\mathcal{L}_nf=0$ then

\begin{multline*}
\mathcal{L}_{m+1}\mathcal{L}_{m+2}\dotsb\mathcal{L}_N\mathcal{R}_{N-1}\mathcal{R}_{N-2}\dotsb \mathcal{R}_nf=q^{-N}(1-q^{N-n})(A_h q^{N+n+h-1}-1)\times\\
\times \mathcal{L}_{m+1}\mathcal{L}_{m+2}\dotsb\mathcal{L}_{N-1}\mathcal{R}_{N-2}\mathcal{R}_{N-2}\dotsb \mathcal{R}_nf\\
=(-1)^{N-m}q^{-(N-m)(N+m+1)/2}(q^{m-n+1};q)_{N-m}(A_h q^{n+m+h};q)_{N-m}\mathcal{R}_{m-1}\mathcal{R}_{m-2}\dotsb \mathcal{R}_nf,
\end{multline*}

\noindent
where the first equality is just \eqref{LRRR} and the second equality follows from the inductive hypothesis.

\end{proof}

\begin{lemma}
If $0\leq n\leq m\leq N$, $f_1\in \text{\rm Ker}\mathcal{L}_m$ and $f_2\in \text{\rm Ker}\mathcal{L}_n$,  then

\begin{multline}\label{RRRRRR}
\langle \mathcal{R}_{N-1}\mathcal{R}_{N-2}\dotsb \mathcal{R}_mf_1,\mathcal{R}_{N-1}\mathcal{R}_{N-2}\dotsb \mathcal{R}_nf_2\rangle_{V_{h,N}}\\
=\delta_{n,m}q^{-(N-n)(N+n+1)/2}(q;q)_{N-n}
(A_hq^{2n+h};q)_{N-n}\langle f_1,f_2\rangle_{V_{h,n}}.
\end{multline}

\end{lemma}

\begin{proof}
It is an immediate consequence of \eqref{invscalarprod} and \eqref{LLLRRR}, noting also that \eqref{invscalarprod} implies that $\text{Ran}\mathcal{R}_{m-1}\bot\text{Ker}\mathcal{L}_m$.
\end{proof}

\begin{corollary}

Suppose that $N\geq n\geq 0$. Then for all $f\in \text{\rm Ker}\mathcal{L}_n$ we have

\begin{equation}\label{multnormRN1Rnh}
\lVert \mathcal{R}_{N-1}\mathcal{R}_{N-2}\dotsb \mathcal{R}_nf \rVert^2_{V_{h,N}}=q^{-(N-n)(N+n+1)/2}(q;q)_{N-n}(A_h q^{2n+h};q)_{N-n}\lVert f\rVert_{V_{h,N}}^2.
\end{equation}

\noindent
In particular, the linear map 

\[
\begin{array}{ccc}
\text{\rm Ker}\mathcal{L}_n&\longrightarrow& V_{h,N}\\
f&\longmapsto&\mathcal{R}_{N-1}\mathcal{R}_{N-2}\dotsb \mathcal{R}_nf
\end{array}
\]

\noindent
is injective.
\end{corollary}

Now we can give the complete spectral analysis of the $q$-Hahn operator $\mathcal{D}_N$.

\begin{theorem}\label{specanmultqHahn}
The following

\begin{equation}\label{multorthdec}
V_{h,N}=\bigoplus_{n=0}^N \mathcal{R}_{N-1}\mathcal{R}_{N-2}\dotsb \mathcal{R}_n\bigl[ \text{\rm Ker}\mathcal{L}_n\bigr]
\end{equation}

\noindent
is the orthogonal decomposition of $V_{h,N}$ into eigenspaces of $\mathcal{D}_N$. The eigenvalue associated to the eigenspace $\mathcal{R}_{N-1}\mathcal{R}_{N-2}\dotsb \mathcal{R}_n\bigl[ \text{\rm Ker}\mathcal{L}_n\bigr]$ 
is equal to 

\[
q^{-n}(1-q^n)(1-A_h q^{n+h-1}).
\]

\end{theorem}

\begin{proof}
For $f\in \text{\rm Ker}\mathcal{L}_n$ we have

\[
\begin{split}
\mathcal{D}_N\bigl[\mathcal{R}_{N-1}\mathcal{R}_{N-2}\dotsb \mathcal{R}_nf\bigr]=
&\mathcal{L}_{N+1}\mathcal{R}_N\mathcal{R}_{N-1}\dotsb \mathcal{R}_nf
+(1-q^{-N-1})(A_h q^{N+h}-1)\mathcal{R}_{N-1}\mathcal{R}_{N-2}\dotsb \mathcal{R}_nf\\
=&q^{-n}(1-q^n)(1-A_h q^{n+h-1})\mathcal{R}_{N-1}\mathcal{R}_{N-2}\dotsb \mathcal{R}_nf,
\end{split}
\]

\noindent
where the first equality follows from \eqref{LN1RN} and the second equality from \eqref{LRRR} (or \eqref{LLLRRR}) (and an elementary calculation). Therefore the $\mathcal{R}_{N-1}\mathcal{R}_{N-2}\dotsb \mathcal{R}_n\bigl[ \text{\rm Ker}\mathcal{L}_n\bigr]$ is an eigenspace of $\mathcal{D}_N$ and the corresponding eigenvalue is $q^{-n}(1-q^n)(1-A_h q^{n+h-1})$. The decomposition \eqref{multorthdec} may be obtained by induction on $N-n$, using the linear algebra identity $V_{h,m}=\text{Ran}\mathcal{R}_{m-1}\bigoplus \text{Ker}\mathcal{L}_m$ and taking into account \eqref{RRRRRR} (the orthogonality also follows from the selfadjointness of $\mathcal{D}_N$).
\end{proof}

\noindent
The preceding Theorem gives the complete spectral analysis of the operator $\mathcal{D}_N$ but in the multidimensional case $h\geq 3$ the eigenspaces are not one-dimensional. In the following sections we will examine the case $h=2$, when the eigenspaces are one-dimensional, and we will show the well known fact that the eigenfunctions are the one variable $q$-Hahn polynomials. In Sections \ref{building} and \ref{treemethod} we will
show that the tree method is a natural way to construct orthogonal bases of eigenfunctions of $\mathcal{D}_N$ in each eigenspace when $h\geq3$. 

%%%%%%%%%%%%%%%%%%%%%%%%%%%%%%%%%%%%%%%%%%%%%%%%%%%%%%%%%%%%%%%%%%%%%%%%%%%%%%%%%%
%%%%%%%%%%%%%%%%%%%%%%%%%%%%%%%%%%%%%%%%%%%%%%%%%%%%%%%%%%%%%%%%%%%%%%%%%%%%%%%%%%%%%
\section{One-dimensional $q$-Hahn polynomials}\label{secqHahn}
%%%%%%%%%%%%%%%%%%%%%%%%%%%%%%%%%%%%%%%%%%%%%%%%%%%%%%%%%%%%%%%%%%%%%%%%%%%%%%%%%%%%%

In this Section we set $h=2$ and we use the following notation: $\alpha=\alpha_1$, $\beta=\alpha_2$, $x=x_1$ and $x_2=N-x_1$. We write a function $f\in V_{2,N}$ in the form $f(x)$ rather than $f(x,N-x)$. 
Now the scalar product \eqref{scalarhmult} has the form

\begin{equation}\label{scalarqHahn}
\langle f_1,f_2\rangle_{V_{2,N}}=q^{N(N+1)/2}\sum_{x=0}^N
\frac{(\alpha q,;q)_x(\beta q;q)_{N-x}}{(q;q)_x(q;q)_{N-x}}(\alpha q)^{N-x}f_1(x)\overline{f_2(x)}
\end{equation}

\noindent
and using the transformation formula $\frac{(\beta q;q)_N(q^{-N};q)_x}{(\beta^{-1}q^{-N};q)_x(q;x)_N}=\frac{\beta^x(\beta q;q)_{N-x}}{(q;q)_{N-x}}$, it is easy to check that it coincides with the usual scalar product for the $q$-Hahn polynomials (see \cite{GaRabook,Ismail,KoSw})) multiplied by $q^{N(N+3)/2}\alpha^N\frac{(\beta q;q)_N}{(q;q)_N}$. Moreover, the $q$-Hahn operator now has the usual form

\[
\mathcal{D}_Nf(x)=B(x)f(x+1)-[B(x)+D(x)]f(x)+D(x)f(x-1)
\]

\noindent
for $f\in V_{h,N}$ and all $x=0,1,2,\dotsc,N$, where 

\[
B(x)=(1-\alpha q^{x+1})(1-q^{x-N})\qquad\quad\text{and}\qquad\quad D(x)=q\alpha(1-q^x)(\beta-q^{x-N-1}).
\]
\noindent
We will need the following transformation formula.

\begin{lemma}\label{RN1Rnh}
Suppose that $0\leq n\leq N$. Then for all $f\in V_{2,N}$ and $0\leq x\leq N$ we have

\begin{equation}\label{RN1Rneqh}
\mathcal{R}_{N-1}\mathcal{R}_{N-2}\dotsb \mathcal{R}_nf(x)=(q;q)_{N-n}\sum_{y=(x-N+n)\vee 0}^{x\wedge n}
\left[\!\!\begin{array}{c}N-x\\n-y\end{array}\!\!\right]_q
\left[\!\!\begin{array}{c}x\\y\end{array}\!\!\right]_q q^{y(y+N-x-n)}q^{-(N-n)(N+n+1)/2}f(y).
\end{equation}

\end{lemma}

\begin{proof}
The proof is by induction on $N$. For $N=n$ it is trivial; it is also easy to check that for $N=n+1$ it coincides with the definition of $\mathcal{R}_n$ (when $h=2$). Assume that \eqref{RN1Rneqh} is true. Then using the elementary identities 
$\left[\!\!\begin{array}{c}m\\k\end{array}\!\!\right]_q=\frac{1-q^{m-k+1}}{1-q^{m+1}}\left[\!\!\begin{array}{c}m+1\\k\end{array}\!\!\right]_q$, 
$(a;q)_{k+1}=(a;q)_k(1-aq^k)$ and the definition of $\mathcal{R}_N$ we get:

\begin{multline*}
\mathcal{R}_{N}\mathcal{R}_{N-1}\dotsb \mathcal{R}_nf(x)=(q;q)_{N-n}\sum_{y=(x-N+n)\vee 0}^{x\wedge n}
\left[\!\!\begin{array}{c}N-x\\n-y\end{array}\!\!\right]_q
\left[\!\!\begin{array}{c}x\\y\end{array}\!\!\right]_q
 q^{y(y+N-x-n)}q^{-(N-n)(N+n+1)/2}\\
 \quad\times q^{x-N-1}(1-q^{-x+N+1})f(y)\\
+(q;q)_{N-n}\sum_{y=(x-N+n-1)\vee 0}^{(x-1)\wedge n}
\left[\!\!\begin{array}{c}N-x+1\\n-y\end{array}\!\!\right]_q
\left[\!\!\begin{array}{c}x-1\\y\end{array}\!\!\right]_q
q^{y(y+N-x-n+1)}q^{-(N-n)(N+n+1)/2}q^{-N-1}
(1-q^x)f(y)\\
=(q;q)_{N-n+1}\sum_{y=(x-N+n-1)\vee 0}^{x\wedge n}
\left[\!\!\begin{array}{c}N-x+1\\n-y\end{array}\!\!\right]_q
\left[\!\!\begin{array}{c}x\\y\end{array}\!\!\right]_q q^{y(y+N-x-n+1)}q^{-(N-n+1)(N+n+2)/2}f(y).
\end{multline*}

\end{proof}

\noindent
Now we are ready to derive the one dimensional $q$-Hahn polynomials as an explicit complete orthogonal system of eigenfunctions for the $q$-difference operator $\mathcal{D}_N$. The results that we give are well known and the method of proof is just a purely analytic version of the group theoretical methods developed by Delsarte \cite{Delsarte}, Dunkl \cite{DuqHahn} and Stanton \cite{Stanton} (a recent account is in Chapter 8 of the monograph \cite{CST3}; see also \cite{MaPa}). Therefore we just sketch the proof. On the other hand, we need to illustrate our methods also in the case $h=2$ because such case form the basis of the recursion procedure in the tree method. We first give a non standard formula and then we derive the classical $_3\phi_2$ expression.\\

\begin{theorem}\label{qHahn}
For $0\leq n\leq N$ and $0\leq x\leq N$ set

\begin{equation}\label{Hahnphi}
Q_n(x;\alpha,\beta,N|q)=\frac{q^{(N-n)(N-n+1)/2}}{(q;q)_{N-n}}\mathcal{R}_{N-1}\mathcal{R}_{N-2}\dotsb \mathcal{R}_n\varphi(x),
\end{equation}  

\noindent
where 

\begin{equation}\label{Qnn}
\varphi(x)\equiv Q_n(x;\alpha,\beta,n|q)=q^{-n^2/2}(q;q)_n(\alpha\beta q^{n+1})^x\frac{(\beta^{-1}q^{-n};q)_x}{(\alpha q;q)_x}.
\end{equation} 

\noindent
Then the $q$-polynomials $Q_n$ satisfy the following $q$-difference equations:

\begin{equation}\label{RN1Qn}
\mathcal{R}_{N-1}Q_n(x;\alpha,\beta,N-1|q)=(q^{-N+n}-1)Q_n(x;\alpha,\beta,N|q),
\end{equation}

\begin{equation}\label{LN1Qn}
\mathcal{L}_{N+1}Q_n(x;\alpha,\beta,N+1|q)=q^{-n}(\alpha\beta q^{N+n+2}-1)Q_n(x;\alpha,\beta,N|q),
\end{equation}

\begin{equation}\label{DNQn}
\mathcal{D}_NQ_n(x;\alpha,\beta,N|q)=q^{-n}(1-q^n)(1-\alpha\beta q^{n+1})Q_n(x;\alpha,\beta,N|q),
\end{equation}

\noindent
have the following special values:

\begin{equation}\label{Qn0}
Q_n(0;\alpha,\beta,N|q)=q^{-n(2N-n)/2}(q^{N-n+1};q)_n,
\end{equation}

\begin{equation}\label{QnN}
Q_n(N;\alpha,\beta,N|q)=q^{-n(2N-n)/2}(q^{N-n+1};q)_n (\alpha\beta q^{n+1})^n\frac{(\beta^{-1}q^{-n};q)_n}{(\alpha q;q)_n},
\end{equation}

\noindent
and satisfy the following orthogonality relations:

\begin{equation}\label{qHahnortrel}
\langle Q_n(\cdot;\alpha,\beta,N|q),Q_m(\cdot;\alpha,\beta,N|q) \rangle_{V_{2,N}}=\delta_{n,m}\frac{(\alpha\beta q^{n+1};q)_{N+1}(q,\beta q;q)_n}{(1-\alpha\beta q^{2n+1})(q;q)_{N-n}(\alpha q;q)_n}\alpha^nq^{[(N-2n)^2+N+2n-2n^2]/2}.
\end{equation}

\noindent
In particular, the polynomials $Q_n$, $n=0,1,\dotsc,N$, form a complete orthogonal system of eigenfunctions for the $q$-difference operator $\mathcal{D}_N$.

\end{theorem}

\begin{proof} 
The function $\varphi(x)$ is the unique solution of the first order equation $\mathcal{L}_n\varphi=0$ satisfying the initial condition $\varphi(0)=q^{-n^2/2}(q;q)_n$. Moreover, using the scalar product \eqref{scalarqHahn} we get:

\[
\begin{split}
\lVert \varphi \rVert^2_{V_{2,n}}= &q^{-n(n-3)/2}[(q;q)_n]^2\alpha^n\sum_{x=0}^n\frac{(\alpha q;q)_x(\beta q;q)_{n-x}}{(q;q)_x(q;q)_{n-x}}(\alpha q)^{-x}(\alpha\beta q^{n+1})^{2x}\left[\frac{(\beta^{-1}q^{-n};q)_x}{(\alpha q;q)_x}\right]^2\\
=&q^{-n(n-3)/2}\alpha^n(q,\beta q;q)_n\sum_{x=0}^n\frac{(\beta^{-1}q^{-n},q^{-n};q)_x}{(q,\alpha q;q)_x}\alpha^x\beta^x q^{x(2n+1)}\\
=&q^{-n(n-3)/2}\alpha^n(q,\beta q;q)_n \;_2\phi_1(q^{-n},\beta^{-1}q^{-n};\alpha q;q,\alpha\beta q^{2n+1})\\
=&q^{-n(n-3)/2}\alpha^n\frac{(q,\beta q,\alpha\beta q^{n+1};q)_n}{(\alpha q;q)_n},
\end{split}
\]

\noindent
where the second equality follows from the transformation formulas $\;(\beta^{-1}q^{-n};q)_x(\beta q;q)_{n-x}=\beta^{-x}(-1)^x$ $\times q^{-x(n-x)-x(x+1)/2}(\beta q;q)_n\;$ and 
$\;\frac{1}{(q;q)_{n-x}}=\frac{(-1)^xq^{x(n-x)+x(x+1)/2}(q^{-n};q)_x}{(q;q)_n}$,
and the fourth equality from the $q$-Vandermonde identity (formula (1.5.2) in \cite{GaRabook}, with $b=\beta^{-1}q^{-n}$ and $c=\alpha q$). \noindent
Finally, using \eqref{multnormRN1Rnh}, \eqref{Hahnphi} and the elementary transformation formulas $\frac{(q;q)_{N-n}}{(q;q)_N}=\frac{(-1)^n}{q^{-n(n-1)/2+Nn}(q^{-N};q)_n}$ and $(\alpha\beta q^{n+1};q)_n(\alpha\beta q^{2n+2};q)_{N-n}=\frac{(\alpha\beta q^{n+1};q)_{N+1}}{1-\alpha\beta q^{2n+1}}$, one can get easily the expression for the norm of $Q_n$ in \eqref{qHahnortrel}.\\

The first order $q$-difference equation \eqref{RN1Qn} follows immediately from \eqref{Hahnphi}, while \eqref{LN1Qn} requires an application of \eqref{LRRR} (or \eqref{LLLRRR}):

\[
\begin{split}
\mathcal{L}_{N+1}Q_n(x;\alpha,\beta,N+1|q)=&
\frac{q^{(N-n+1)(N-n+2)/2}}{(q;q)_{N-n+1}}\mathcal{L}_{N+1}\mathcal{R}_N\mathcal{R}_{N-1}\dotsb \mathcal{R}_n\varphi(x)\\
=&\frac{q^{(N-n+1)(N-n+2)/2}}{(q;q)_{N-n+1}}\cdot q^{-N-1}(1-q^{N-n+1})(\alpha\beta q^{N+n+2}-1)\mathcal{R}_{N-1}\dotsb \mathcal{R}_n\varphi(x)\\
=&q^{-n}(\alpha\beta q^{N+n+2}-1)Q_n(x;\alpha,\beta,N|q).
\end{split}
\]

\noindent
Finally, \eqref{DNQn} may be deduced from \eqref{RN1Qn}, \eqref{LN1Qn} and \eqref{RN1LN} but it is also a particular case of the spectral analysis in Theorem \ref{specanmultqHahn}.

\end{proof}

Now we derive the classical $_3 \phi_2$-expression for $Q_n$.

\begin{proposition}
The $q$-Hahn polynomials have the following classical $\;_3\phi_2$ expression:

\[
Q_n(x;\alpha,\beta,N|q)=q^{-nN+n^2/2}\frac{(q;q)_N}{(q;q)_{N-n}}
\;_3\phi_2
\left[\begin{array}{c}
q^{-n},\alpha\beta q^{n+1},q^{-x}\\
\alpha q,q^{-N}\end{array};q,q
\right]
\]
\end{proposition}
\begin{proof}
From \eqref{RN1Rneqh}, \eqref{Hahnphi} and \eqref{Qnn} we get easily

\[
Q_n(x;\alpha,\beta,N|q)=q^{-nN+n^2/2}(q^{N-x-n+1};q)_n\;_3\phi_2
\left[\begin{array}{c}
q^{-n},q^{-x},\beta^{-1}q^{-n}\\
\alpha q,q^{N-x-n+1}\end{array};q,\alpha\beta q^{N+n+2}
\right]
\]

\noindent
and then the Proposition follows from an application of formula (3.2.5) in \cite{GaRabook}, with $a=q^{-x}$, $b={\beta^{-1}q^{-n}}$, $d=\alpha q$ and $e=q^{N-x-n+1}$.
\end{proof}

\begin{remark}
{\rm 
Usually, in the literature the $q$-Hahn polynomials are denoted by $Q_n(q^{-x};\alpha,\beta,N|q)$ and are equal precisely to $\;_3\phi_2
\left[\begin{array}{c}
q^{-n},\alpha\beta q^{n+1},q^{-x}\\
\alpha q,q^{-N}\end{array};q,q
\right]$. We use a different normalization (and also a slightly different notation) in order to get simpler expression for the action of the lowering and raising operators and for the connection coefficients between two-dimensional $q$-Hahn polynomials.
Note also that the norm in \eqref{qHahnortrel} coincides with the usual norm for the $q$-Hahn polynomials multiplied by $q^{N(N+3)/2}\alpha^N\frac{(\beta q;q)_N}{(q;q)_N}\left[q^{-nN+n^2/2)}(q^{N-n+1};q)_n\right]^2$. In \cite{GaRa} another different notation (and normalization) is used for the $q$-Hahn polynomials, namely they are denoted by the symbol $h_n(x;a,b,N,;q)$ (formula (3.16) in \cite{GaRa}); since we must compare our results with those in \cite{GaRa}, we state explicitly the relation with our notation:

\begin{equation}\label{GaRaqHahn}
Q_n(x;\alpha,\beta,N|q)=\frac{(-1)^{N-n}q^{n/2}}{(\alpha q;q)_n}h_n(x;\alpha,\beta,N;q).
\end{equation}
}
\end{remark}

In our recursive definition of the multidimensional $q$-Hahn polynomials we will use the following $x_1,x_2$-notation for the one-dimensional case:

\begin{equation}\label{qHahnx1x2}
f_{n,x_1+x_2}(x_1,x_2)=Q_n(x_1;\alpha_1,\alpha_2,x_1+x_2|q). 
\end{equation}

\noindent 
In particular, \eqref{RN1Qn} and \eqref{LN1Qn} have the following explicit expressions: 

\begin{equation}\label{RqHahn}
\begin{split}
q^{-x_1-x_2}(1-q^{x_1})&Q_n(x_1-1;\alpha_1,\alpha_2,x_1+x_2-1|q)+q^{-x_2}(1-q^{x_2}) Q_n(x_1;\alpha_1,\alpha_2,x_1+x_2-1|q)\\
=&(q^{-x_1-x_2+n}-1)Q_n(x_1;\alpha_1,\alpha_2,x_1+x_2|q),
\end{split}
\end{equation}

\begin{equation}\label{LqHahn}
\begin{split}
(\alpha_1 q^{x_1+1}-1)&Q_n(x_1+1;\alpha_1,\alpha_2,x_1+x_2+1|q)+\alpha_1q^{x_1+1}(\alpha_2 q^{x_2+1}-1)Q_n(x_1;\alpha_1,\alpha_2,x_1+x_2+1|q)\\
=&q^{-n}(\alpha_1\alpha_2q^{x_1+x_2+n+2}-1)Q_n(x_1;\alpha_1,\alpha_2,x_1+x_2|q),
\end{split}
\end{equation}

\noindent
Note that \eqref{RqHahn} and \eqref{LqHahn} coincide respectively with i) and ii) in \cite{DuqHahn}, Proposition 2.4. \\

We end this section with a simple identity that will be used in Section \ref{bidqHahn}.

\begin{lemma}
We have:

\begin{equation}\label{QnVandermonde}
\sum_{x=0}^j(-1)^xq^{x(x-1)/2}\frac{Q_n(x;\alpha,\beta,n|q)}{(q;q)_x(q;q)_{j-x}}=\frac{(\alpha\beta q^{n+1};q)_j}{(\alpha q,q;q)_j}q^{-n^2/2}(q;q)_n.
\end{equation}
\end{lemma}

\begin{proof}
Using \eqref{Qnn} and the transformation formula $(q;q)_{j-x}=(-1)^xq^{-jx+x(x-1)/2)}\frac{(q;q)_j}{(q^{-j};q)_x}$ we get

\[
\sum_{x=0}^j(-1)^xq^{x(x-1)/2}\frac{Q_n(x;\alpha,\beta,n|q)}{(q;q)_x(q;q)_{j-x}}=q^{-n^2/2}\frac{(q;q)_n}{(q;q)_j}\cdot\;_2\phi_1
\left[\begin{array}{c}
q^{-j},\beta^{-1} q^{-n}\\
\alpha q \end{array};q,\alpha\beta q^{n+j+1}
\right].
\]

\noindent
Then \eqref{QnVandermonde} follows from an application of the $q$-Vandermonde identity (1.5.2) in \cite{GaRabook} (with $n=j$, $b=\beta^{-1}q^{-n}$ and $c=\alpha q$).
\end{proof}

%%%%%%%%%%%%%%%%%%%%%%%%%%%%%%%%%%%%%%%%%%%%%%%%%%%%%%%%%%%%%%%%%%%%%%%%%%%%%%%%%%
%%%%%%%%%%%%%%%%%%%%%%%%%%%%%%%%%%%%%%%%%%%%%%%%%%%%%%%%%%%%%%%%%%%%%%%%%%%%%%%%%%%%%
\section{Building the tree method on the factorization method}\label{building}
%%%%%%%%%%%%%%%%%%%%%%%%%%%%%%%%%%%%%%%%%%%%%%%%%%%%%%%%%%%%%%%%%%%%%%%%%%%%%%%%%%%%%

In this section we give some preliminary but fundamental results involving raising and lowering operators that we will subsequently use to develop the tree method. First of all, we need to introduce (and explain progressively) some particular notation. We fix two positive integers $t,h$, with $1\leq t\leq h-1$, and for any $h$-parts composition ${\bf x}=(x_1,x_2,\dotsc,x_h)$ we set ${\bf x'}=(x_1,x_2,\dotsc,x_t)$ and ${\bf x}''=(x_{t+1},x_{t+1},\dotsc,x_h)$. We have the following decomposition (where $\coprod$ means disjoint union):

\begin{equation}\label{dechN}
[h;N]=\coprod_{M=0}^N\Bigl\{[t;M]\times[h-t;N-M]\Bigr\}
\end{equation}

\noindent
obtained simply by writing ${\bf x}=({\bf x'},{\bf x''})$. From \eqref{dechN} we deduce the following decomposition of $V_{h,N}$:

\begin{equation}\label{decLhN}
V_{h,N}=\bigoplus_{M=0}^N \Bigl[ V_{t,M}\otimes V_{h-t,N-M}\Bigr].
\end{equation}

\noindent
If $f'_M\in V_{t,M}$ and $f''_{N-M}\in V_{h-t,N-M}$ then the tensor product $f'_M\otimes f''_{N-M}$ is defined by setting $(f'_M\otimes f''_{N-M})({\bf x'},{\bf x''})=f'_M({\bf x'})\cdot f''_{N-M}({\bf x''})$ for all ${\bf x'}\in [t;M]$ and ${\bf x''}\in [h-t,N-M]$. A more suggestive form is: $f({\bf x})=f'_{X_t}({\bf x'})f''_{X_h-X_t}({\bf x''})$, replacing $M$ with $X_t$ and $N-M$ with $X_h-X_t$. This more intrinsic notation may be used also for $f$: we can write $f_{X_h}({\bf x})$, allowing $N\equiv X_h$ to vary. Note also that $f'_{X_t}$ and $f_{X_h-X_t}$ are determined up to a multiplicative constant that depends only on $(X_h,X_t-X_h)$. Therefore we may consider functions $f_{X_h}$ of the form 

\begin{equation}\label{basicfact}
f_{X_h}({\bf x})=\psi(X_t,X_h-X_t)f'_{X_t}({\bf x'})f''_{X_h-X_t}({\bf x''}),
\end{equation}

\noindent
where $f_{X_t}$ and $f_{X_h-X_t}$ are defined respectively on $V_{t,X_t}$ and $V_{h-t,X_h-X_t}$, for certain values of $X_t$ and $X_h-X_t$, and $\psi$ will be a function of the two (numerical) variables $(X_h,X_t-X_h)$.\\

We denote by $\mathcal{L}_M', \mathcal{R}_M'$ (respectively $\mathcal{L}_{N-M}'', \mathcal{R}_{N-M}''$) the lowering and raising operators defined on $V_{t,M}$ with parameters $\alpha_1,\alpha_2,\dotsc,\alpha_t$ (respectively defined on $V_{h-t,N-M}$ with parameters $\alpha_{t+1},\alpha_{t+2},\dotsc,\alpha_h$). For instance, if $f'\in V_{t,M}$ then

\[
\mathcal{L}'_Mf'({\bf x'})=\sum_{j=1}^tA_{j-1}q^{j-1+X_{j-1}}(\alpha_jq^{x_j+1}-1)f'(\epsilon_j^+({\bf x})),\qquad\text{for all }\quad{\bf x'}\in[t;M-1],
\]

\noindent
while if $f''\in V_{h-t,N-M}$ then

\[
\mathcal{L}''_{N-M}f''({\bf x''})=\sum_{j={t+1}}^hA_{j-1}A_t^{-1}q^{j-t-1+X_{j-1}-X_t}(\alpha_jq^{x_j+1}-1)f''(\epsilon_j^+({\bf x})),\qquad\text{for all }\quad{\bf x''}\in[h-t;N-M-1].
\]

\noindent
In the degenerate case $t=1$ we take $f'_{x_1}$ constant and we define $\mathcal{L}'$ and $\mathcal{R}'$ by setting $\mathcal{L}'_{x_1}=(\alpha_1q^{x_1+1}-1)I$ and $\mathcal{R}'_{x_1}=q^{-x_1}(1-q^{x_1})I$, where $I$ is the identity. Similarly, if $t=h-1$ we take $f''_{x_h}$ constant, $\mathcal{L}''_{x_h}=(\alpha_hq^{x_h+1}-1)I$ and $\mathcal{R}''_{x_h}=q^{-x_h}(1-q^{x_h})I$.

\begin{proposition}
Suppose that $f_{X_h}({\bf x})$ is as in \eqref{basicfact}. Then we have

\begin{equation}\label{Lfpsi}
\begin{split}
\mathcal{L}_{X_h+1}f_{X_h+1}({\bf x})=&\psi(X_t+1,X_h-X_t)\mathcal{L}'_{X_t+1}f'_{X_t+1}({\bf x'})\cdot f''_{X_h-X_t}({\bf x''})\\
&+A_tq^{t+X_t}\psi(X_t,X_h-X_t+1) f'_{X_t}({\bf x'})\cdot\mathcal{L}''_{X_h-X_t+1}f'_{X_h-X_t+1}({\bf x''})
\end{split}
\end{equation}

\noindent
and

\begin{equation}\label{Rfpsi}
\begin{split}
\mathcal{R}_{X_h-1}f_{X_h-1}({\bf x})=&q^{-X_h+X_t} \psi(X_t-1,X_h-X_t)\mathcal{R}'_{X_t-1}f'_{X_t-1}({\bf x'})\cdot f''_{X_h-X_t}({\bf x''})\\
&+\psi(X_t,X_h-X_t-1) f'_{X_t}({\bf x'})\cdot\mathcal{R}''_{X_h-X_t-1}f'_{X_h-X_t-1}({\bf x''}),
\end{split}
\end{equation}

\noindent
for all ${\bf x}=({\bf x'},{\bf x''})$.
\end{proposition}

\begin{proof}
From the definitions of $\mathcal{L}_N$ and $\mathcal{R}_N$ we get immediately

\begin{multline*}
\mathcal{L}_{X_h+1}f_{X_h+1}({\bf x})=\sum_{j=1}^tA_{j-1}q^{j-1+X_{j-1}}(\alpha_jq^{x_j+1}-1)\psi(X_t+1,X_h-X_t)f_{X_t+1}'(\epsilon_j^+({\bf x'})) f''_{X_h-X_t}({\bf x''})\\
+A_tq^{t+X_t}\sum_{j=t+1}^hA_{j-1}A_t^{-1}q^{(j-t-1)+X_{j-1}-X_t}(\alpha_jq^{x_j+1}-1)
\psi(X_t,X_h-X_t+1) f'_{X_t}({\bf x'}) f_{X_h-X_t+1}''(\epsilon_j^+({\bf x''}))\\
=\psi(X_t+1,X_h-X_t)\mathcal{L}'_{X_t+1}f'_{X_t+1}({\bf x'})\cdot f''_{X_h-X_t}({\bf x''})\qquad\qquad\qquad\\
+A_tq^{t+X_t}\psi(X_t,X_h-X_t+1) f'_{X_t}({\bf x'})\cdot\mathcal{L}''_{X_h-X_t+1}f''_{X_h-X_t+1}({\bf x''})\qquad\qquad\qquad\qquad
\end{multline*}

\noindent
and

\[
\begin{split}
\mathcal{R}_{X_h-1}f_{X_h-1}({\bf x})=&
q^{-X_h+X_t}\sum_{i=1}^tq^{X_{i-1}-X_t}(1-q^{x_i})\psi(X_t-1,X_h-X_t) f'_{X_t-1}(\epsilon_i^-({\bf x'})) f''_{X_h-X_t}({\bf x''})\\
&+\sum_{i=t+1}^hq^{X_{i-1}-X_h}(1-q^{x_i})\psi(X_t,X_h-X_t-1)f'_{X_t}({\bf x'}) f''_{X_h-X_t-1}(\epsilon_i^-({\bf x''}))\\
=&q^{-X_h+X_t}\psi(X_t-1,X_h-X_t)\mathcal{R}'_{X_t-1}f'_{X_t-1}({\bf x'})\cdot f''_{X_h-X_t}({\bf x''})\\
&+\psi(X_t,X_h-X_t-1) f'_{X_t}({\bf x'})\cdot\mathcal{R}''_{X_h-X_t-1}f''_{X_h-X_t-1}({\bf x''}).
\end{split}
\]

\end{proof}

\begin{proposition}\label{scalarprodfxgx}
Suppose that $f_{X_h}$ is as in \eqref{basicfact} and that similarly $g_{X_h}$ has the form $g_{X_h}({\bf x})=\phi(X_t,X_h-X_t)g'_{X_t}({\bf x'})g''_{X_h-X_t}({\bf x''})$. Then we have:

\begin{equation}\label{fNgN}
\begin{split}
\langle f_N, g_N \rangle_{V_{h,N}}=&\sum_{M=0}^N(A_t q^{M+t})^{N-M}\langle f'_M, g'_M \rangle_{V_{t,M}}\langle f''_{N-M}, g''_{N-M} \rangle_{V_{h-t,N-M}}\times\\
&\times\psi(M,N-M)\overline{\phi(M,N-M)}.
\end{split}
\end{equation}

\end{proposition}

\begin{proof}
Indeed, since $\sum_{{\bf x}\in [h;N]}=\sum_{M=0}^N\sum_{{\bf x'}\in [t;M]}\sum_{{\bf x''}\in [h-t;N-M]}$, from the expression of the scalar product \eqref{scalarhmult} we get

\[
\begin{split}
\langle f_N,&g_N \rangle_{V_{h,N}}=q^{N(N+1)/2}\sum_{M=0}^N\Biggl\{\sum_{{\bf x'}\in [t;M]}f'_M({\bf x'}) \overline{g'_M({\bf x'})}\prod_{i=1}^t\Biggl[\frac{(q\alpha_i;q)_{x_i}}{(q;q)_{x_i}}\left(\alpha_i q\right)^{X_t-X_i}\Biggr]\Biggr\}\\
&\times\Biggl[\prod_{i=1}^t\left(\alpha_iq\right)^{X_h-X_t}\Biggr]\cdot
\Biggl\{\sum_{{\bf x''}\in [h-t;N-M]}
f''_{N-M}({\bf x''})\overline{g'_{N-M}({\bf x''})}
\prod_{i=t+1}^h\Biggl[\frac{(q\alpha_i;q)_{x_i}}{(q;q)_{x_i}}\left(\alpha_i q\right)^{X_h-X_i}\Biggr]\Biggr\}\\
&\times\psi(M,N-M)\overline{\phi(M,N-M)}.
\end{split}
\]

\noindent
The expressions in curly brackets coincides respectively with 

\[
q^{-M(M+1)/2}\langle f'_M, g'_M \rangle_{V_{t,M}} \qquad\text{and}\qquad
q^{-(N-M)(N-M+1)/2}\langle f''_{N-M}, g''_{N-M} \rangle_{V_{h-t,N-M}},
\]

\noindent while $q^{N(N+1)/2}q^{-M(M+1)/2}q^{-(N-M)(N-M+1)/2}\prod_{i=1}^t\left(\alpha_iq\right)^{X_h-X_t}=(\alpha_1\alpha_2\dotsb \alpha_t q^{M+t})^{N-M}$.
\end{proof}

Now we make a precise choice of the function $\psi(X_t,X_h-X_t)$ and give some recursive results for the action of the lowering and raising operators and for the scalar product. We give a purely analytic version of Theorem 4.19 in Dunkl's paper \cite{DuqHahn}; see also Section 2.3 of our paper \cite{ScarabottimultHahn} for the case $q=1$.

\begin{theorem}\label{theoremrecursion}
Suppose that $i,j,n$ are nonnegative integers satisfying $i+j\leq n$ and
in \eqref{basicfact} take $\psi(X_t,X_h-X_t)=q^{-jX_t}\psi_{n-i-j}(X_t,X_h-X_t)$, where 

\begin{equation}\label{psiQn}
\psi_{n-i-j}(X_t,X_h-X_t)=Q_{n-i-j}(X_t-i;A_tq^{t+2i-1},A_hA_t^{-1}q^{h-t+2j-1},X_h-i-j|q).
\end{equation}

\noindent
If $t=1$ we always take $i=0$ while if $t=h-1$ we always take $j=0$.

\begin{enumerate}
\item
If the functions $f'_{X_t}$ and $f''_{X_h-X_t}$ satisfy the first order $q$-difference identities

\[
\mathcal{L}'_{X_t+1}f'_{X_t+1}=q^{-i}(A_tq^{t+X_t+i}-1)f'_{X_t}
\]

\noindent
and

\[
\mathcal{L}''_{X_h-X_t+1}f''_{X_h-X_t+1}=q^{-j}(A_hA_t^{-1}q^{h-t+X_h-X_t+j}-1)f''_{X_h-X_t}
\]

\noindent
then $f_{X_h}$ satisfies

\begin{equation}\label{Lfx}
\mathcal{L}_{X_h+1}f_{X_h+1}=q^{-n}(A_hq^{h+X_h+n}-1)f_{X_h}.
\end{equation}

\item
If the functions $f'_{X_t}$ and $f''_{X_h-X_t}$ satisfy the first order $q$-difference equations

\[
\mathcal{R}'_{X_t-1}f'_{X_t-1}=(q^{-X_t+i}-1)f'_{X_t}\qquad \text{and}\qquad \mathcal{R}''_{X_h-X_t-1}f''_{X_h-X_t-1}=(q^{-X_h+X_t+j}-1)f''_{X_h-X_t}
\]

\noindent
then $f_{X_h}$ satisfies

\begin{equation}\label{Rfx}
\mathcal{R}_{X_h-1}f_{X_h-1}=(q^{-X_h+n}-1)f_{X_h}.
\end{equation}
\end{enumerate}

\end{theorem}

\begin{proof}
\begin{enumerate}
\item
We have

\[
\begin{split}
\mathcal{L}_{X_h+1}f_{X_h+1}({\bf x})=&f'_{X_t}({\bf x'})f_{X_h-X_t}''({\bf x''})q^{-i-j}q^{-jX_t}\Bigl[(A_t q^{t+X_t+i}-1)\psi_{n-i-j}(X_t+1,X_h-X_t)\\
&+A_tq^{t+X_t+i}(A_hA_t^{-1}q^{h-t+X_h-X_t+j}-1)\psi_{n-i-j}(X_t,X_h-X_t+1)
\Bigr]\\
=&q^{-n}(A_hq^{h+X_h+n}-1)f_{X_h}({\bf x}),
\end{split}
\]

\noindent
where we have applied \eqref{Lfpsi} in the first identity and \eqref{LqHahn} in the second identity.

\item
Similarly, applying \eqref{Rfpsi} and \eqref{RqHahn} we have

\[
\begin{split}
\mathcal{R}_{X_h-1}f_{X_h-1}({\bf x})=&f'_{X_t}({\bf x'})f_{X_h-X_t}''({\bf x''})q^{-jX_t}\Bigl[q^{-X_h+i+j}(1- q^{X_t-i})\psi_{n-i-j}(X_t-1,X_h-X_t)\\
&+q^{-X_h+X_t+j}(1-q^{X_h-X_t-j})\psi_{n-i-j}(X_t,X_h-X_t-1)
\Bigr]\\
=&(q^{-X_h+n}-1)f_{X_h}({\bf x}),
\end{split}
\]

\end{enumerate}

\noindent
If $t=1$ or $t=h-1$ the finite difference identities \eqref{Lfx} and \eqref{Rfx} are verified in virtue of the definitions of the operators itself. 

\end{proof}

\begin{corollary}\label{corollaryrecursion}
If $f_{X_h}$ satisfies \eqref{Lfx} and \eqref{Rfx} then

\[
\mathcal{D}_{X_h}f_{X_h}=q^{-n}(1-q^n)(1-A_h q^{h+n-1})f_{X_h}.
\] 
\end{corollary}

\begin{proof}
It is an immediate consequence of \eqref{RN1LN} (or of \eqref{LN1RN}).
\end{proof}

\begin{theorem}\label{scalarprodnm}
Fix $N$ and two nonnegative integers $i,i$ satisfying $i+j\leq N$. Suppose that $f'_M$ and $f''_{N-M}$ are defined for $i\leq M\leq N-j$ and that

\begin{equation}\label{normfM}
\lVert f'_M\rVert^2_{V_{t,M}}=\Gamma'\cdot\frac{(A_t q^{t+2i};q)_{M-i}}{(q;q)_{M-i}}q^{[(M-2i)^2+M+2i-2i^2]/2},
\end{equation}

\begin{equation}\label{normfNM}
\lVert f''_{N-M}\rVert^2_{V_{h-t,N-M}}=\Gamma''\cdot\frac{(A_hA_t^{-1} q^{h-t+2j};q)_{N-M-j}}{(q;q)_{N-M-j}}q^{[(N-M-2j)^2+N-M+2j-2j^2]/2},
\end{equation}

\noindent
where the constant $\Gamma'$ and $\Gamma''$ do not depend on $M$. Let $n,m$ be two nonnegative integers satisfying the conditions $i+j\leq n\leq N$ and $i+j\leq m\leq N$. Let $\psi_{n-i-j}$ be as in \eqref{psiQn} (with $x'+x''=N$) and set

\[
f({\bf x})=q^{-jX_t}\psi_{n-i-j}(X_t,X_h-X_t)f_{X_t}({\bf x'})f_{X_h-X_t}({\bf x''})\equiv q^{-jM}\psi_{n-i-j}(M,N-M)f_M({\bf x'})f_{N-M}({\bf x''}),
\]

\[
g({\bf x})=q^{-jX_t}\psi_{m-i-j}(X_t,X_h-X_t)f_{X_t}({\bf x'})f_{X_h-X_t}({\bf x''})\equiv q^{-jM}\psi_{m-i-j}(M,N-M)f_M({\bf x'})f_{N-M}({\bf x''}),
\]

\noindent
for all ${\bf x}\in [h;N]$ such that $i\leq X_t\leq N-j$; for all the other values of ${\bf x}$ set $f({\bf x})=g({\bf x})=0$. Then we have:

\begin{equation}\label{fgVhN}
\langle f,g\rangle_{V_{h,N}}=\delta_{n,m}\cdot \Gamma\cdot\frac{(A_h q^{h+2n};q)_{N-n}}{(q;q)_{N-n}}q^{[(N-2n)^2+N+2n-2n^2]/2},
\end{equation}

\noindent
where

\begin{equation}\label{Xi}
\Gamma=\Gamma'\cdot \Gamma''\cdot \frac{(q,A_h q^{h+n+i+j-1},A_hA_t^{-1} q^{h-t+2j};q)_{n-i-j}}{(A_t q^{t+2i};q)_{n-i-j}}(A_t q^{t+2i})^{n-i}q^{-2ij-(n-i-j)},
\end{equation}

\noindent
(and therefore $\Gamma$ does not depend on $N$).

\end{theorem}

\begin{proof}
We have

\[
\begin{split}
\langle f,g\rangle_{V_{h,N}}=&\Gamma'\cdot \Gamma''\cdot q^{i(i-2N+1)+j(j-2N+1+t+2i)+\frac{N(N+1)}{2}}A_t^j\sum_{M=i}^{N-j}\psi_{n-i-j}(M,N-M)\\
&\times\psi_{m-i-j}(M,N-M)\frac{(A_tq^{t+2i};q)_{M-i}(A_hA_t^{-1}q^{h-t+2j};q)_{N-M-j}}{(q;q)_{M-i}(q;q)_{N-M-j}}(A_tq^{t+2i})^{N-M-j}\\
=&\delta_{n,m}\cdot \Gamma'\cdot \Gamma''\cdot(A_t q^{t+2i})^{n-i}q^{-2ij-(n-i-j)}\frac{(q,A_hA_t^{-1} q^{h-t+2j};q)_{n-i-j}}{(A_t q^{t+2i};q)_{n-i-j}}\\
&\times q^{[(N-2n)^2+N+2n-2n^2]/2} \frac{(A_h q^{h+n+i+j-1};q)_{N-i-j+1}}{(q;q)_{N-n}(1-A_h q^{h+2n-1})}\\
=& \delta_{n,m}\cdot\Gamma\cdot\frac{(A_h q^{2h+n};q)_{N-n}}{(q;q)_{N-n}}q^{[(N-2n)^2+N+2n-2n^2]/2},
\end{split}
\]

\noindent
where in the first equality we have used \eqref{fNgN} and the hypothesis \eqref{normfM}, \eqref{normfNM}, in the second equality we have used the orthogonality relations \eqref{qHahnortrel} (and we have also rearranged the powers of $q$ and $A_t$) and in the final equality we have used \eqref{Xi} and the identity

\begin{equation}\label{identitynorm}
\frac{(A_h q^{h+n+i+j-1};q)_{N-i-j+1}}{1-A_h q^{h+2n-1}}=(A_h q^{h+n+i+j-1};q)_{n-i-j}(A_h q^{h+2n};q)_{N-n}.
\end{equation}

\end{proof}

%%%%%%%%%%%%%%%%%%%%%%%%%%%%%%%%%%%%%%%%%%%%%%%%%%%%%%%%%%%%%%%%%%%%%%%%%%%%%%%%%%
%%%%%%%%%%%%%%%%%%%%%%%%%%%%%%%%%%%%%%%%%%%%%%%%%%%%%%%%%%%%%%%%%%%%%%%%%%%%%%%%%%%%%
\section{The tree method for multidimensional $q$-Hahn polynomials}\label{treemethod}
%%%%%%%%%%%%%%%%%%%%%%%%%%%%%%%%%%%%%%%%%%%%%%%%%%%%%%%%%%%%%%%%%%%%%%%%%%%%%%%%%%%%%

A {\em rooted binary tree} $\mathcal{T}$ is a tree with a distinguished vertex $V$ ({\em the root}) of degree 2 and all the remaining vertices of degree 3 or 1. The vertices of degree 1 are called {\em leaves}, all the other vertices (including the root) are called {\em internal vertices}, or {\em branch points}. The $l$-th {\em level} of a tree $\mathcal{T}$, denoted by $\mathcal{T}_l$, is formed by the vertices at distance $l$ from the root. The {\em height} of $\mathcal{T}$ is the greatest $L$ such that there exists a vertex in $\mathcal{T}$ at distance $L$ from the root. If $U\in\mathcal{T}_l$ is an internal vertex, then there exist exactly two vertices $X,Y\in\mathcal{T}_{l+1}$ connected with $V$; they are called the {\em sons} of $U$, while $U$ is the {\em father} of $X$ and $Y$. We think of $\mathcal{T}$ as a {\em planar} tree, and therefore $U$ has a {\em left} son and a {\em right} son.  In the figure below, $X$ is the left son and $Y$ is the right son.

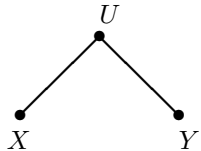
\begin{picture}(400,80)
\put(200,50){\circle*{4}}
\put(170,20){\circle*{4}}
\put(230,20){\circle*{4}}
\put(200,55){$U$}
\put(165,7){$X$}
\put(230,7){$Y$}

\thicklines
\put(200,50){\line(-1,-1){30}}
\put(200,50){\line(1,-1){30}}

\end{picture}

For a tree $\mathcal{T}$, we denote by $\mathcal{T}'$ (resp. $\mathcal{T}''$) the subtrees formed by the left (resp. right) descendants of the root. 
We always denote by $V$ the root of $\mathcal{T}$ and by $W$ and $Z$ respectively its left and right son. Then $W$ is the root of $\mathcal{T}'$ and $Z$ is the root of $\mathcal{T}''$. There is a basic recursive procedure that we will use many times: if something has been proved/defined for $\mathcal{T}'$ and $\mathcal{T}''$, then we can use this fact to prove/define the same thing for $\mathcal{T}$ (see also our previous paper \cite{ScarabottimultHahn}).\\

\begin{proposition}\label{intvert}
In a rooted binary tree with $h$ leaves, the number of internal points is equal to $h-1$.
\end{proposition}

\begin{proof}
The Proposition is obvious for $h=2$. We may prove the general case by induction, observing that if $\mathcal{T}'$ has $t$ leaves and $t-1$ internal vertices and $\mathcal{T}''$ has $h-t$ leaves and $h-t-1$ internal vertices then $\mathcal{T}$ has $t+(h-t)=h$ leaves and $(t-1)+(h-t-1)+1=h-1$ internal vertices. 
\end{proof}

Now suppose that $\mathcal{T}$ has $h$ leaves and fix a set of {\em parameters} $\alpha_1,\alpha_2,\dotsc,\alpha_h$ satisfying the conditions $0<\alpha_i<q^{-1}$, $i=1,2,\dotsc,h$, or the conditions $\alpha_i>q^{-N}$, $i=1,2,\dotsc,h$. Let $x_1,x_2,\dotsc, x_h$ be a set of {\em variables}. The associated {\em parameters labeling} is constructed in the following way: we label the root $V$ with $(\alpha_1,\alpha_2,\dotsc,\alpha_h)$; then if $t$ is the number of leaves of $\mathcal{T}'$, we label $W$ with $(\alpha_1,\alpha_2,\dotsc,\alpha_t)$ and $Z$ with $(\alpha_{t+1},\alpha_{t+2},\dotsc,\alpha_h)$:

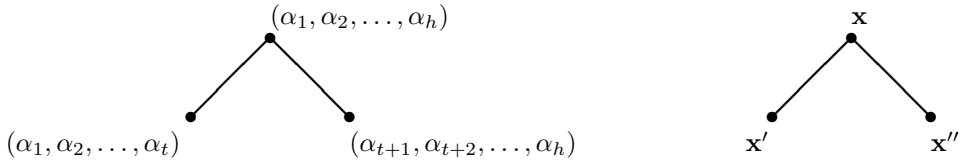
\begin{picture}(400,80)
\put(100,50){\circle*{4}}
\put(70,20){\circle*{4}}
\put(130,20){\circle*{4}}
\put(100,55){$(\alpha_1,\alpha_2,\dotsc,\alpha_h)$}
\put(0,7){$(\alpha_1,\alpha_2,\dotsc,\alpha_t)$}
\put(130,7){$(\alpha_{t+1},\alpha_{t+2},\dotsc,\alpha_h)$}

\thicklines
\put(100,50){\line(-1,-1){30}}
\put(100,50){\line(1,-1){30}}

\put(320,50){\circle*{4}}
\put(290,20){\circle*{4}}
\put(350,20){\circle*{4}}
\put(320,55){${\bf x}$}
\put(280,7){${\bf x'}$}
\put(350,7){${\bf x''}$}

\thicklines
\put(320,50){\line(-1,-1){30}}
\put(320,50){\line(1,-1){30}}
\end{picture}

\noindent
and then we can iterate this procedure. Similarly, we can construct an an associated {\em variables labeling}: we label $V,W,Z$ respectively with ${\bf x},{\bf x'},{\bf x''}$ and we proceed recursively. This way every internal vertex is labeled with a sequence of consecutive $\alpha$'s (or $x$'s) while every leaf is labeled with a single $\alpha$ (or $x$). The following is an example with $h=4$.

\begin{picture}(450,150)
\thicklines
%%%%%%%%%%%%%%% zero level
%%%%%%%%%%%%%%% first tree
\put(100,120){\circle*{4}}
%%%%%% labelings
\put(105,120){$(\alpha_1,\alpha_2,\alpha_3,\alpha_4)$}

%%%%%%%%%%%%%%second tree
\put(350,120){\circle*{4}}
%%%%%% labelings
\put(355,120){$(x_1,x_2,x_3,x_4)$}

%%%%%%%%%%%% first level 

%%%%%%%%%%%%%%%  first tree
\put(160,60){\circle*{4}}
\put(40,60){\circle*{4}}

\put(100,120){\line(-1,-1){90}}
\put(100,120){\line(1,-1){90}}

%%%%%%%%%%%% labelings
\put(165,60){$(\alpha_3,\alpha_4)$}
\put(0,60){$(\alpha_1,\alpha_2)$}

%%%%%%%%%%%%% second tree
\put(410,60){\circle*{4}}
\put(290,60){\circle*{4}}

\put(350,120){\line(-1,-1){90}}
\put(350,120){\line(1,-1){90}}

%%%%%%%%%%%% labelings
\put(415,60){$(x_3,x_4)$}
\put(252,60){$(x_1,x_2)$}

%%%%%%%%%%%%%%%%%%% second level

%%%%%%%%%%%%%%%%%%% first tree
\put(130,30){\circle*{4}}
\put(190,30){\circle*{4}}
\put(70,30){\circle*{4}}
\put(10,30){\circle*{4}}

\put(40,60){\line(1,-1){30}}
\put(160,60){\line(-1,-1){30}}

%%%%%%%%%%%%%%%%%% labelings
\put(130,20){$\alpha_3$}
\put(190,20){$\alpha_4$}
\put(65,20){$\alpha_2$}
\put(5,20){$\alpha_1$}

%%%%%%%%%%%%%%%%%% second tree

\put(380,30){\circle*{4}}
\put(440,30){\circle*{4}}
\put(320,30){\circle*{4}}
\put(260,30){\circle*{4}}

\put(290,60){\line(1,-1){30}}
\put(410,60){\line(-1,-1){30}}
 %%%%%%%%%%%%%%%%%% labelings
\put(380,20){$x_3$}
\put(440,20){$x_4$}
\put(315,20){$x_2$}
\put(255,20){$x_1$}

\put(200,0){{\bf Figure 1}}
\end{picture}

\qquad\\
\qquad\\

A {\em coefficients labeling} $\text{c}$ for $\mathcal{T}$ is defined by assigning a nonnegative integer to each internal vertex, 0 to each leaf. If $U\in\mathcal{T}$, $\text{c}(U)$ is the coefficient associated to $U$.
We denote by $\text{CL}(\mathcal{T},n)$ the set of all coefficients labellings $c$ of $\mathcal{T}$ such that $\sum_{U\in\mathcal{T}}\text{c}(U)=n$. The following is an example of a coefficients labeling. 

\begin{picture}(450,150)
\thicklines
%%%%%%%%%%%%%%% zero level
\put(100,120){\circle*{4}}
%%%%%% labelings
\put(105,120){$n-i-j$}

%%%%%%%%%%%% first level 

\put(160,60){\circle*{4}}
\put(40,60){\circle*{4}}

\put(100,120){\line(-1,-1){90}}
\put(100,120){\line(1,-1){90}}

%%%%%%%%%%%% labelings
\put(167,60){$j$}
\put(30,60){$i$}

%%%%%%%%%%%%%%%%%%% second level

\put(130,30){\circle*{4}}
\put(190,30){\circle*{4}}
\put(70,30){\circle*{4}}
\put(10,30){\circle*{4}}

\put(40,60){\line(1,-1){30}}
\put(160,60){\line(-1,-1){30}}

%%%%%%%%%%%%%%%%%% labelings
\put(130,20){$0$}
\put(190,20){$0$}
\put(65,20){$0$}
\put(5,20){$0$}

\put(80,0){{\bf Figure 2}}
\end{picture}

\quad\\

The following proposition is an immediate consequence of Proposition \ref{intvert}. 

\begin{proposition}\label{dimVhN}
The cardinality of the set $\coprod_{n=0}^N\text{CL}(\mathcal{T},n)$ is equal to the cardinality of $[h;N]$.
\end{proposition}

If $U$ is a vertex of $\mathcal{T}$ and $(\alpha_{l+1},\alpha_{l+2},\dotsc,\alpha_{m})$  (resp. $(x_{l+1},x_{l+2},\dotsc,x_m)$) is the parameters label of $U$ (resp. its variables label), we set $\text{p}(U)=\alpha_{l+1}\alpha_{l+2}\dotsb\alpha_{m}q^{m-l}$ (resp. $\text{v}(U)=x_{l+1}+x_{l+2}+\dotsb+x_m$).
If $X$ and $Y$ are respectively the left and the right son of $U$, we set $\text{lp}(U)=\text{p}(X)$, $\text{rp}(U)=\text{p}(Y)$, $\text{lv}(U)=\text{v}(X)$ and $\text{rv}(U)=\text{v}(Y)$. We also denote by $\text{lcs}(U)$ (resp. $\text{rcs}(U)$) the sum of all the coefficients of the left (resp. right) descendants of $U$ (lcs={\em left coefficients sum}, while rcs={\em right coefficients sum}) and we set $\text{cs}(U)=\text{c}(U)+\text{lcs}(U)+\text{rcs}(U)$ (that is, $\text{cs}(U)$ is the sum of the coefficients of all the vertices of the subtree rooted at $U$).\\

Now we are in position to define the multidimensional $q$-Hahn polynomials associated to a rooted tree $\mathcal{T}$. Suppose that $\text{c}\in \text{CL}(\mathcal{T},n)$ and that $i$ (resp. $j$) is the sum of the labels of the vertices in $\mathcal{T}'$ (resp. $\mathcal{T}''$). Then $n-i-j$ is the label of the root, and we denote by $\text{c}'$ (resp. $\text{c}''$) the coefficients labeling of $\mathcal{T}'$ (resp. $\mathcal{T}''$). Then the multidimensional $q$-Hahn polynomials associated to $\mathcal{T}(n)$ is defined in the following recursive way:

\begin{equation}\label{multqHahn}
\begin{split}
Q_{\text{c}}({\bf x};\alpha_1,\dotsc,\alpha_h,x|q)=&q^{-jX_t}Q_{n-i-j}(X_t-i;A_tq^{t+2i-1},A_hA_t^{-1}q^{h-t+2j-1},X_h-i-j|q)\\
&\times Q_{\text{c}'}({\bf x'};\alpha_1,\dotsc,\alpha_t,X_t|q)\cdot Q_{\text{c}''}({\bf x''};\alpha_{t+1},\dotsc,\alpha_h,X_h-X_t|q).\\
\end{split}
\end{equation}

\noindent
If $t=1$ we set $Q_{\text{c}'}({\bf x'};\alpha_1,\dotsc,\alpha_t,x_1|q)=1$, while if $t=h-1$ we set $Q_{\text{c}''}({\bf x''};\alpha_{t+1},\dotsc,\alpha_h,X_h-X_{h-1}|q)=1$. Moreover, if $h=2$ then necessarily $i=j=0$, $t=1$ and \eqref{multqHahn} coincides with \eqref{qHahnx1x2}. Note also that we must have $X_h\geq n$, $X_t\geq i$ and $X_h-X_t\geq j$; applying recursively these conditions we find that \eqref{multqHahn} is defined for those ${\bf x}\in [h;N]$ such that: 

\begin{equation}\label{existcond}
\text{v}(U)\geq \text{cs}(U),\quad\qquad \text{for all}\quad U\in\mathcal{T}.
\end{equation}

\noindent
For the values of ${\bf x}$ that do not satisfy the conditions \eqref{existcond}, we set $Q_{\text{c}}({\bf x};\alpha_1,\dotsc,\alpha_h,x|q)=0$. Finally, we define a real valued function $\Gamma$ by setting, for every internal vertex $U\in\mathcal{T}$,

\begin{equation}\label{XiU}
\Gamma(U)=\frac{(q,\text{p}(U)q^{\text{cs}(U)+\text{lcs}(U)+\text{rcs}(U)-1},\text{rp}(U)q^{2\text{rcs}(U)};q)_{\text{c}(U)}}{(\text{lp}(U)q^{2\text{lcs}(U)};q)_{\text{c}(U)}}[\text{lp}(U)q^{2\text{lcs}(U)}]^{\text{c}(U)+\text{rcs}(U)}
q^{-2\text{lcs}(U)\text{rcs}(U)-\text{c}(U)},
\end{equation}

\noindent
and $\Gamma(U)=1$ if $U$ is a leaf. Clearly, $\Gamma(U)$ is modeled on \eqref{Xi}.\\

Now we can state the first fundamental result of the present paper, in which we give the main properties of the multidimensional $q$-Hahn polynomials.

\begin{theorem}\label{theomultqHahn}
\begin{enumerate}
\item
The set $\{Q_{\text{\rm c}}({\bf \cdot};\alpha_1,\dotsc,\alpha_h,N|q):0\leq n\leq N,\text{\rm c}\in \text{\rm CL}(\mathcal{T},n)\}$ is an orthogonal basis for $V_{h,N}$.

\item\label{theomultqHahneigen}
For $0\leq n\leq N$, the set $\{Q_{\text{\rm c}}({\bf \cdot};\alpha_1,\dotsc,\alpha_h,N|q):\text{\rm c}\in \text{\rm CL}(\mathcal{T},n)\}$ is an orthogonal basis for 

\[
\mathcal{R}_{N-1}\mathcal{R}_{N-2}\dotsb \mathcal{R}_n\bigl[\text{\rm Ker}\mathcal{L}_n\bigr],
\]

\noindent
that is for the eigenspace of $\mathcal{D}_N$ corresponding to the eigenvalue $q^{-n}(1-q^n)(1-A_hq^{h+n-1})$.

\item
We have 

\begin{equation}\label{normQc}
\lVert Q_{\text{\rm c}}({\bf \cdot};\alpha_1,\alpha_2,\dotsc,\alpha_h,N|q)\rVert_{V_{h,N}}^2=
\frac{(A_h q^{h+2n};q)_{N-n}}{(q;q)_{N-n}}q^{[(N-2n)^2+N+2n-2n^2]/2}\prod_{U\in\mathcal{T}}\Gamma(U).
\end{equation}

\item
The polynomials $Q_{\text{\rm c}}({\bf \cdot};\alpha_1,\dotsc,\alpha_h,N|q)$, $\text{\rm c}\in \text{\rm CL}(\mathcal{T},n)$, satisfy the following recurrence relation:

\begin{equation}\label{Qcmultform}
Q_{\text{\rm c}}({\bf x};\alpha_1,\dotsc,\alpha_h,N|q)=\frac{q^{(N-n)(N-n+1)/2}}{(q;q)_{N-n}}\mathcal{R}_{N-1}\mathcal{R}_{N-2}\dotsb \mathcal{R}_nQ_{\text{\rm c}}({\bf x};\alpha_1,\dotsc,\alpha_h,n|q).
\end{equation}

\end{enumerate}
\end{theorem}

\begin{proof} These results follows form our recursive definition of $Q_\text{c}$, taking the function $f_{n,x_1+x_2}$ in \eqref{qHahnx1x2} as the basis of the induction. First of all, note that for $h=2$ the expression \eqref{normQc} becomes

\[
\alpha_1^n\frac{(q,\alpha_1\alpha_2 q^{n+1},\alpha_2q;q)_n(\alpha_1\alpha_2q^{2n+2};q)_{N-n}}{(\alpha_1q;q)_n(q;q)_{N-n}}q^{[(N-2n)^2+N+2n-2n^2]/2},
\]

\noindent
and this coincides with $\lVert f_{n,N}\rVert^2_{V_{2,N}}$, which is given by \eqref{qHahnortrel} (one has just to apply \eqref{identitynorm} for $h=2$ and $i=j=0$). \\

Now suppose that $\text{c}_1, \text{c}_2$ are two different labellings of $\mathcal{T}$.
Our recursive definition yields $\langle Q_{\text{c}_1},Q_{\text{c}_2}\rangle_{V_{h,N}}=0$: we can apply Proposition \ref{scalarprodfxgx} when $\text{c}_1'\neq\text{c}_2'$ or $\text{c}_1''\neq\text{c}_2''$, while when $\text{c}_1'=\text{c}_2'$ and $\text{c}_1''=\text{c}_2''$ but clearly $\text{c}_1(V)\neq\text{c}_2(V)$, we can invoke \eqref{fgVhN} in Theorem \ref{scalarprodnm}. The norm of $Q_\text{c}$ may be computed using recursively \eqref{fgVhN} and \eqref{Xi}, and this lead to \eqref{XiU} and \eqref{normQc} as the final result. The orthogonal system $\{Q_{\text{c}}({\bf \cdot};\alpha_1,\dotsc,\alpha_h,N|q):\text{c}\in \text{CL}(\mathcal{T},n),0\leq n\leq N\}$ is complete in $V_{h,N}$ simply because $\text{dim}V_{h,N}$ is equal to the cardinality of $\coprod_{n=0}^N\text{CL}(\mathcal{T},n)$ (Proposition \ref{dimVhN}).\\ 

Similarly, 2. follows from a recursive application of Theorem \ref{theoremrecursion}, taking again the functions \eqref{qHahnx1x2} as the basis of the recursion, and applying Corollary \ref{corollaryrecursion}. Finally, \eqref{Qcmultform} follows from a repeated application of \eqref{Rfx}.
 
\end{proof}

It is possible to characterize the members of our orthogonal basis for $V_{h,N}$ as the common eigenfunctions of a set of $q$-difference operators. We need to introduce other notation and definitions. Suppose again that $U$ is an internal vertex of $\mathcal{T}$ and that $(\alpha_{l+1},\alpha_{l+2},\dotsc,\alpha_m)$ (resp. $(x_{l+1},x_{l+1},\dotsc,x_m)$) is its parameter label (resp. variable label). We associate to $U$ the multidimensional $q$-Hahn operator $\mathcal{D}_U$ defined as in \eqref{multqHahnoperator} but acting on functions of the variables $x_{l+1},x_{l+2},\dotsc,x_m$, with parameters $\alpha_{l+1},\alpha_{l+2},\dotsc,\alpha_m$ ans $x_{l+1}+x_{l+2}+\dotsb+x_m$ in place of $N$. Then our recursive definition of $Q_\text{c}$ together with \ref{theomultqHahneigen}. in Theorem \ref{theomultqHahn} give immediately the following proposition.

\begin{proposition}
The polynomial $Q_\text{c}$ is an eigenfunction of $\mathcal{D}_U$ and the corresponding eigenvalue is equal to $\lambda_{\text{\rm c},U}=q^{-\text{\rm cs}(U)}(1-q^{\text{\rm cs}(U)})(1-\text{\rm p}(U)q^{\text{\rm cs}(U)-1})$. Moreover the set $\{\lambda_{\text{\rm c},U}:U \text{ is an internal vertex of }\mathcal{T}\}$ characterizes $Q_\text{\rm c}$.
\end{proposition}

\begin{example}\label{3dimqHahn}
{\rm
Consider the tree in Figure 1 and denote by $\text{c}$ the labeling in Figure 2. The associated three-dimensional $q$-Hahn polynomial is given by:

\[
\begin{split}
Q_{\text{c}}(x_1,x_2,&x_3,x_4;\alpha_1,\alpha_2,\alpha_3,\alpha_4,x_1+x_2+x_3+x_4|q)=q^{-j(x_1+x_2)}Q_{n-i-j}
(x_1+x_2-i;\alpha_1\alpha_2q^{2i+1},\\
&\alpha_3\alpha_4q^{2j+1},x_1+x_2+x_3+x_4-i-j|q)\cdot
Q_i(x_1;\alpha_1,\alpha_2,x_1+x_2|q)\cdot Q_j(x_3;\alpha_3,\alpha_4,x_3+x_4|q),
\end{split}
\]

\noindent
for $(x_1,x_2,x_3,x_4)\in [4,N]$.
}
\end{example}

\begin{example}\label{bidqHahn1}{\rm
Consider now the following tree with the parameters, variables and coefficients labellings depicted below.

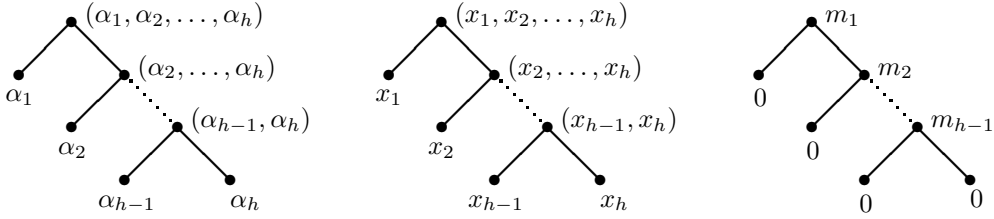
\begin{picture}(400,150)
%%%%%%%%%%%%%%% first tree

%%%%%%%%%%%%%%% zero level
\put(20,120){\circle*{4}}
%%%%%% labelings
\put(25,120){$(\alpha_1,\alpha_2,\dotsc,\alpha_h)$}

%%%%%%%%%%%% first level 
\put(0,100){\circle*{4}}
\put(40,100){\circle*{4}}

\thicklines
\put(20,120){\line(-1,-1){20}}
\put(20,120){\line(1,-1){20}}
\dottedline{4}(40,100)(60,80)

%%%%%%%%%% labelings
\put(-5,90){$\alpha_1$}
\put(45,100){$(\alpha_2,\dotsc,\alpha_h)$}

%%%%%%%%%%%%%%%  second level
\put(20,80){\circle*{4}}
\put(60,80){\circle*{4}}

\put(60,80){\line(-1,-1){20}}

%%%%%%%%%%%% labelings
\put(65,80){$(\alpha_{h-1},\alpha_h)$}
\put(15,70){$\alpha_2$}

%%%%%%%%%%%%%%%%%%% third level
\put(80,60){\circle*{4}}
\put(40,60){\circle*{4}}

\put(60,80){\line(1,-1){20}}
\put(40,100){\line(-1,-1){20}}
%%%%%%%%%%%%%%%% labelings
\put(30,50){$\alpha_{h-1}$}
\put(80,50){$\alpha_h$}

%%%%%%%%%%%%%%% second tree

%%%%%%%%%%%%%%% zero level
\put(160,120){\circle*{4}}
%%%%%% labelings
\put(165,120){$(x_1,x_2,\dotsc,x_h)$}

%%%%%%%%%%%% first level 
\put(140,100){\circle*{4}}
\put(180,100){\circle*{4}}

\thicklines
\put(160,120){\line(-1,-1){20}}
\put(160,120){\line(1,-1){20}}
\dottedline{4}(180,100)(200,80)

%%%%%%%%%% labelings
\put(135,90){$x_1$}
\put(185,100){$(x_2,\dotsc,x_h)$}

%%%%%%%%%%%%%%%  second level
\put(160,80){\circle*{4}}
\put(200,80){\circle*{4}}

\put(200,80){\line(-1,-1){20}}

%%%%%%%%%%%% labelings
\put(205,80){$(x_{h-1},x_h)$}
\put(155,70){$x_2$}

%%%%%%%%%%%%%%%%%%% third level
\put(220,60){\circle*{4}}
\put(180,60){\circle*{4}}

\put(200,80){\line(1,-1){20}}
\put(180,100){\line(-1,-1){20}}
%%%%%%%%%%%%%%%% labelings
\put(170,50){$x_{h-1}$}
\put(220,50){$x_h$}

%%%%%%%%%%%%%%% third tree

%%%%%%%%%%%%%%% zero level
\put(300,120){\circle*{4}}
%%%%%% labelings
\put(305,120){$m_1$}

%%%%%%%%%%%% first level 
\put(280,100){\circle*{4}}
\put(320,100){\circle*{4}}

\thicklines
\put(300,120){\line(-1,-1){20}}
\put(300,120){\line(1,-1){20}}
\dottedline{4}(320,100)(340,80)

%%%%%%%%%% labelings
\put(278,88){$0$}
\put(325,100){$m_2$}

%%%%%%%%%%%%%%%  second level
\put(300,80){\circle*{4}}
\put(340,80){\circle*{4}}

\put(340,80){\line(-1,-1){20}}

%%%%%%%%%%%% labelings
\put(345,80){$m_{h-1}$}
\put(298,68){$0$}

%%%%%%%%%%%%%%%%%%% third level
\put(360,60){\circle*{4}}
\put(320,60){\circle*{4}}

\put(340,80){\line(1,-1){20}}
\put(320,100){\line(-1,-1){20}}
%%%%%%%%%%%%%%%% labelings
\put(318,48){$0$}
\put(360,50){$0$}

\end{picture}

\noindent
Set $j_k=m_{k+1}+m_{k+2}+\dotsb+m_{h-1}$, for $k=0,1,\dotsc,h-2$. Then the associated $(h-1)$-dimensional $q$-Hahn polynomial, that we denote by $\xi_{m_1,m_2,\dotsc,m_{h-1}}({\bf x};\alpha_1,\dotsc,\alpha_h,X_h|q)$, is given by the following formula:

\[
\begin{split}
\xi_{m_1,m_2,\dotsc,m_{h-1}}({\bf x};\alpha_1,&\dotsc,\alpha_h,X_h|q)=q^{-j_1x_1}Q_{m_1}(x_1;\alpha_1,A_hA_1^{-1}q^{h+2j_1-2},X_h-j_1|q)\\
&\times q^{-x_2j_2}Q_{m_2}(x_2;\alpha_2,A_hA_2^{-1}q^{h+2j_2-3},X_h-X_1-j_2|q)\times\dotsb\\ 
\dotsb&\times q^{-x_{h-2}j_{h-2}}Q_{m_{h-2}}(x_{h-2};\alpha_{h-2},\alpha_{h-1}\alpha_hq^{2j_{h-2}+1},x_{h-2}+x_{h-1}+x_h-j_{h-2}|q)\\	
&\times Q_{m_{h-1}}(x_{h-1};\alpha_{h-1},\alpha_h,x_{h-1}+x_h|q).
\end{split}
\]

\noindent
For this polynomial the conditions \eqref{existcond} become simply $\quad x_k+x_{k+1}+\dotsb +x_h\geq j_{k-1},\quad k=1,2,\dotsb,h-1$, and their
norm may be obtained applying \eqref{normQc}:

\begin{equation}\label{normxi}
\begin{split}
\lVert \xi_{m_1,m_2,\dotsc,m_{h-1}}\rVert^2_{V_{h,N}}=&\frac{(A_hq^{h+2n};q)_{N-n}}{(q;q)_{N-n}}q^{[(N-2n)^2+N+2n-2n^2]/2}\\
&\times\prod_{k=1}^{h-1}\frac{(q,A_hA_{k-1}^{-1}q^{h-k+j_{k-1}+j_k},A_hA_k^{-1}q^{h-k+2j_k};q)_{m_k}}{(\alpha_k q;q)_{m_k}}(\alpha_k q)^{j_{k-1}}q^{-m_k}
\end{split}
\end{equation}

}
\end{example}

\begin{example}\label{bidqHahn2}{\rm
Consider now the following tree with the parameters, variables and coefficients labellings depicted below.

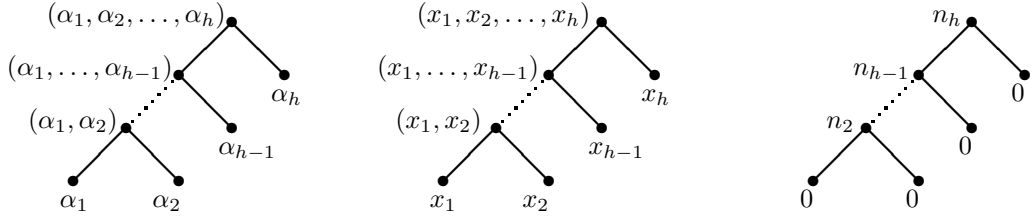
\begin{picture}(400,150)
%%%%%%%%%%%%%%% first tree

%%%%%%%%%%%%%%% zero level
\put(70,120){\circle*{4}}
%%%%%% labelings
\put(0,120){$(\alpha_1,\alpha_2,\dotsc,\alpha_h)$}

%%%%%%%%%%%% first level 
\put(90,100){\circle*{4}}
\put(50,100){\circle*{4}}

\thicklines
\put(70,120){\line(-1,-1){20}}
\put(70,120){\line(1,-1){20}}
\dottedline{4}(50,100)(30,80)

%%%%%%%%%% labelings
\put(85,90){$\alpha_h$}
\put(-15,100){$(\alpha_1,\dotsc,\alpha_{h-1})$}

%%%%%%%%%%%%%%%  second level
\put(30,80){\circle*{4}}
\put(70,80){\circle*{4}}

\put(30,80){\line(-1,-1){20}}

%%%%%%%%%%%% labelings
\put(-7,80){$(\alpha_1,\alpha_2)$}
\put(65,70){$\alpha_{h-1}$}

%%%%%%%%%%%%%%%%%%% third level
\put(10,60){\circle*{4}}
\put(50,60){\circle*{4}}

\put(30,80){\line(1,-1){20}}
\put(50,100){\line(1,-1){20}}
%%%%%%%%%%%%%%%% labelings
\put(40,50){$\alpha_2$}
\put(5,50){$\alpha_1$}

%%%%%%%%%%%%%%% second tree

%%%%%%%%%%%%%%% zero level
\put(210,120){\circle*{4}}
%%%%%% labelings
\put(140,120){$(x_1,x_2,\dotsc,x_h)$}

%%%%%%%%%%%% first level 
\put(230,100){\circle*{4}}
\put(190,100){\circle*{4}}

\thicklines
\put(210,120){\line(-1,-1){20}}
\put(210,120){\line(1,-1){20}}
\dottedline{4}(190,100)(170,80)

%%%%%%%%%% labelings
\put(225,90){$x_h$}
\put(125,100){$(x_1,\dotsc,x_{h-1})$}

%%%%%%%%%%%%%%%  second level
\put(170,80){\circle*{4}}
\put(210,80){\circle*{4}}

\put(170,80){\line(-1,-1){20}}

%%%%%%%%%%%% labelings
\put(132,80){$(x_1,x_2)$}
\put(205,70){$x_{h-1}$}

%%%%%%%%%%%%%%%%%%% third level
\put(150,60){\circle*{4}}
\put(190,60){\circle*{4}}

\put(170,80){\line(1,-1){20}}
\put(190,100){\line(1,-1){20}}
%%%%%%%%%%%%%%%% labelings
\put(180,50){$x_2$}
\put(145,50){$x_1$}

%%%%%%%%%%%%%%% third tree
%%%%%%%%%%%%%%% zero level
\put(350,120){\circle*{4}}
%%%%%% labelings
\put(335,120){$n_h$}

%%%%%%%%%%%% first level 
\put(370,100){\circle*{4}}
\put(330,100){\circle*{4}}

\thicklines
\put(350,120){\line(-1,-1){20}}
\put(350,120){\line(1,-1){20}}
\dottedline{4}(330,100)(310,80)

%%%%%%%%%% labelings
\put(365,90){$0$}
\put(305,100){$n_{h-1}$}

%%%%%%%%%%%%%%%  second level
\put(310,80){\circle*{4}}
\put(350,80){\circle*{4}}

\put(310,80){\line(-1,-1){20}}

%%%%%%%%%%%% labelings
\put(295,80){$n_2$}
\put(345,70){$0$}

%%%%%%%%%%%%%%%%%%% third level
\put(290,60){\circle*{4}}
\put(330,60){\circle*{4}}

\put(310,80){\line(1,-1){20}}
\put(330,100){\line(1,-1){20}}
%%%%%%%%%%%%%%%% labelings
\put(325,50){$0$}
\put(285,50){$0$}

\end{picture}

\noindent
Set $i_k=n_2+n_3+\dotsb+n_{k-1}$, for $k=3,4,\dotsc,h,h+1$. Then the associated $(h-1)$-dimensional $q$-Hahn polynomial, that we denote by $\theta_{n_2,n_3,\dotsc,n_h}({\bf x};\alpha_1,\dotsc,\alpha_h,X_h|q)$, is given by the following formula:

\[
\begin{split}
\theta_{n_2,n_3,\dotsc,n_h}({\bf x};\alpha_1,\dotsc,\alpha_h,X_h|q)=&Q_{n_h}(X_{h-1}-i_h;A_{h-1}q^{h+2i_h-2},\alpha_h,X_h-i_h|q)\\
&\times Q_{n_{h-1}}(X_{h-2}-i_{h-1};A_{h-2}q^{h+2i_{h-1}-3},\alpha_{h-1},X_{h-1}-i_{h-1}|q)\times\dotsb\\ 
\dotsb&\times Q_{n_3}(x_1+x_2-i_3;\alpha_1\alpha_2q^{2i_3+1},\alpha_3,x_1+x_2+x_3-i_3|q)\\	
&\times Q_{n_2}(x_1;\alpha_1,\alpha_2,x_1+x_2|q).
\end{split}
\]

\noindent
For this polynomial the conditions \eqref{existcond} become simply $\quad x_1+x_2+\dotsb+x_{k-1}\geq i_k,\quad k=3,4,\dotsb,h+1$, and again their norm polynomials may be obtained applying \eqref{normQc}:

\begin{equation}\label{normtheta}
\begin{split}
\lVert \theta_{n_2,n_3,\dotsc,n_h}\rVert^2_{V_{h,N}}=&\frac{(A_hq^{h+2n};q)_{N-n}}{(q;q)_{N-n}}q^{[(N-2n)^2+N+2n-2n^2]/2}\\
&\times\prod_{k=2}^h\frac{(q,A_kq^{k+i_{k-1}+i_k-1},\alpha_kq;q)_{n_k}}{(A_{k-1}q^{k+2i_k-1};q)_{n_k}}((A_{k-1}q^{k+2i_k-1})^{n_k+i_k}q^{-n_k}.
\end{split}
\end{equation}

\noindent
The polynomials obtained in this example coincide with those of Gasper and Rahman in \cite{GaRa}, formula (3.15), modulo a different normalization; see \eqref{GaRaqHahn} in the present paper.
}
\end{example}

%%%%%%%%%%%%%%%%%%%%%%%%%%%%%%%%%%%%%%%%%%%%%%%%%%%%%%%%%%%%%%%%%%%%%%%%%%%%%%%%%%
%%%%%%%%%%%%%%%%%%%%%%%%%%%%%%%%%%%%%%%%%%%%%%%%%%%%%%%%%%%%%%%%%%%%%%%%%%%%%%%%%%%%%
\section{$q$-Racah polynomials as connection coefficients between two-dimensional $q$-Hahn polynomials: Dunkl's method}\label{bidqHahn}
%%%%%%%%%%%%%%%%%%%%%%%%%%%%%%%%%%%%%%%%%%%%%%%%%%%%%%%%%%%%%%%%%%%%%%%%%%%%%%%%%%%%%

The present section is a translation in our setting of the results in sections 3. and 4. of Dunkl's paper \cite{Du6}, with some minor modifications in the methods of proof. We set $h=3$ and a function $f\in V_{3,N}$ will be written in the form $f(x_1,x_2)$ (we omit $x_3\equiv N-x_1-x_2$). The domain of definition of a function $f\in V_{3,N}$ is then the set of points with integer coordinates in the triangle of vertices $(0,N),(0,0),(N,0)$ in the $x_1,x_2$-plane (Figure 3a). \\

Now suppose that $f\in V_{3,N}$ and $\mathcal{L}_Nf=0$. The explicit form of this equation is:

\[
(\alpha_1 q^{x_1+1}-1)f(x_1+1,x_2)+\alpha_1 q^{x_1+1}(\alpha_2 q^{x_2+1}-1)f(x_1,x_2+1)+\alpha_1\alpha_2 q^{x_1+x_2+2}(\alpha_3 q^{N-x_1-x_2}-1)f(x_1,x_2)=0.
\]

\noindent
Therefore $f(x_1,x_2)$ is determined by $f(x_1,x_2-1)$ and $f(x_1+1,x_2-1)$. This may be used recursively to prove that each function in $V_{3,N}\cap\text{Ker}\mathcal{L}_N$ is determined by its values at the points $(0,0),(1,0),\dotsc,(N,0)$; in particular, the value of $f(x_1,x_2)$ depends only on the values of $f$ at $(x_1,0), (x_1+1,0),\dotsc,(x_1+x_2,0)$ (the domain of dependence of $(x_1,x_2)$; Figure 3b). Therefore the value $f(k,0)$ has the set $\{(x_1,x_2):0\leq x_1\leq k, k-x_1\leq x_2\leq N-x_1\}$ as its domain of influence (Figure 3c).\\

\begin{picture}(450,140)
\thicklines

%%%%%%%%%%%%%%% first triangle
\put(0,20){\circle*{3}}
\put(100,20){\circle*{3}}
\put(0,120){\circle*{3}}
\put(0,20){\line(1,0){100}}
\put(0,20){\line(0,1){100}}
\put(100,20){\line(-1,1){100}}
\put(90,10){$(N,0)$}
\put(-10,10){$(0,0)$}
\put(-10,125){$(0,N)$}

\put(25,-5){{\bf Figure 3a}}

%%%%%%%%%%%%%%% second triangle
\put(170,20){\circle*{3}}
\put(230,20){\circle*{3}}
\put(170,80){\circle*{3}}
\put(170,60){\circle*{3}}
\put(190,60){\circle*{3}}
\dashline{3}(190,20)(230,20)
\put(170,80){\line(0,-1){20}}
\put(170,80){\line(1,-1){20}}
\dashline{3}(170,60)(170,20)
\dashline{3}(190,60)(230,20)
\put(210,10){$(x_1+x_2,0)$}
\put(160,10){$(x_1,0)$}
\put(150,85){$(x_1,x_2)$}
\put(115,57){$(x_1,x_2-1)$}
\put(195,57){$(x_1+1,x_2-1)$}

\put(175,-5){{\bf Figure 3b}}

%%%%%%%%%%%%%%% third triangle
\put(300,20){\circle*{3}}
\put(400,20){\circle*{3}}
\put(300,120){\circle*{3}}
\put(360,20){\circle*{3}}
\put(300,20){\line(0,1){100}}
\put(300,20){\line(1,0){100}}
\put(400,20){\line(-1,1){100}}
\put(360,20){\line(-1,1){60}}
\put(360,20){\line(0,1){40}}

\thinlines
\put(355,25){\line(0,1){40}}
\put(350,30){\line(0,1){40}}
\put(345,35){\line(0,1){40}}
\put(340,40){\line(0,1){40}}
\put(335,45){\line(0,1){40}}
\put(330,50){\line(0,1){40}}
\put(325,55){\line(0,1){40}}
\put(320,60){\line(0,1){40}}
\put(315,65){\line(0,1){40}}
\put(310,70){\line(0,1){40}}
\put(305,75){\line(0,1){40}}

\put(390,10){$(N,0)$}
\put(290,10){$(0,0)$}
\put(290,125){$(0,N)$}
\put(350,10){$(k,0)$}

\put(325,-5){{\bf Figure 3c}}

\end{picture}
\quad\\

Now we translate in our setting Proposition 2.5 in \cite{Du6}, which gives an explicit representation for functions in $V_{3,N}\cap\text{Ker}\mathcal{L}_N$. 

\begin{proposition}
For $k=0,1,2,\dotsc, N$ set

\[
\begin{split}
f_{N,k}(x_1,x_2)=&\left[\!\!\begin{array}{c}x_2\\k-x_1\end{array}\!\!\right]_q\frac{(\alpha_1q^{x_1+1};q)_{k-x_1}(\alpha_3q^{N-x_1-x_2+1};q)_{x_1+x_2-k}}{(\alpha_2q;q)_{x_2}}\times\\
&\times\alpha_1^{x_1-k}\alpha_2^{x_1+x_2-k}q^{(x_1-k)(x_1+x_2+1)+x_2(x_2+1)/2}.
\end{split}
\]

\noindent
Then every function $f\in V_{3,N}\cap\text{\rm Ker}\mathcal{L}_N$ may be represented in the form

\begin{equation}\label{Duformula}
f(x_1,x_2)=\sum_{k=x_1}^{x_1+x_2}f(k,0)f_{N,k}(x_1,x_2).
\end{equation}

\end{proposition} 

\begin{proof}
For $k=0,1,\dotsc N$ we have:

\[
f_{N,k}\in V_{3,N}\cap\text{\rm Ker}\mathcal{L}_N,
\]

\begin{equation}\label{fNjjk}
f_{N,k}(k,0)=1,\qquad f_{N,k}(h,0)=0 \quad\text{ for }\quad h\neq k.
\end{equation}

\noindent
Hence the right hand side of \eqref{Duformula} belongs to $V_{3,N}\cap\text{\rm Ker}\mathcal{L}_N$ and is equal to $f$ on each point $(k,0)$, $k=0,1,2,\dotsc,N$. Therefore it coincides with $f$. 

\end{proof}

\noindent
Note also that the support of $f_{N,k}$ is precisely the domain of influence of $(k,0)$ (Figure 3c) and that this is a consequence of \eqref{fNjjk}.\\

From Theorem \ref{theomultqHahn}, Examples \ref{bidqHahn1} and \ref{bidqHahn2} we know that the sets $\xi_{n-j,j}$, $j=0,1,2,\dotsc,n$ and $\theta_{n-i,i}$, $i=0,1,2,\dotsc,n$ are two different bases for $\mathcal{R}_{N-1}\mathcal{R}_{N-2}\dotsb \mathcal{R}_n\bigl[\text{\rm Ker}\mathcal{L}_n\bigr]$.
The main goal of this section is to find the connection coefficients between these two bases. As in \cite{Du6}, it suffices to find these coefficients in the case $N=n$, since we can invoke 4. in Theorem \ref{theomultqHahn} to translate the results for $N=n$ to the case $N>n$. We introduce the following notation for the functions $\xi$'s and $\theta$'s in the case $N=n$:

\[
\widetilde{\xi}_{n-j,j}(x_1,x_2)=\xi_{n-j,j}(x_1,x_2,n-x_1-x_2;\alpha_1,\alpha_2\alpha_3|q),
\]

\[
\begin{split}
\widetilde{\theta}_{n-i,i}(x_1,x_2)=&\theta_{n-i,i}(x_1,x_2,n-x_1-x_2;\alpha_1,\alpha_2\alpha_3|q)\\
\equiv& Q_{n-i}(x_1+x_2-i;\alpha_1\alpha_2 q^{2i+1},\alpha_3,n-i|q)\cdot Q_i(x_1;\alpha_1,\alpha_2,x_1+x_2|q).
\end{split}
\]

\noindent
The functions $\xi_{n-j,j}$, $\theta_{n-i,i}$, $\widetilde{\xi}_{n-j,j}$ and $\widetilde{\theta}_{n-i,i}$ correspond respectively to the functions $\widehat{\psi}_{rk}$, $\widehat{\phi}_{rm}$, $\psi_{rk}$ and $\phi_{rm}$ in \cite{Du6}; this may be proved using formulas (2.4) and (2.9) in \cite{Du6} (see also Section 2 in \cite{DuqHahn}). \\

Now we give our version of Theorem 3.1 in \cite{Du6}, with a minor simplification in the proof (in Dunkl's notation, we take $y=r-m$).

\begin{theorem}
Suppose that $f\in\text{\rm Ker}\mathcal{L}_n$ and that $f=\sum_{i=0}^n a_i\widetilde{\theta}_{n-i,i}$. Then we have

\begin{equation}\label{ai}
a_i=\frac{(-\alpha_2)^i(\alpha_1q;q)_iq^{(3i^2+i+n^2)/2-ni}}{(q;q)_{n-i}(\alpha_1\alpha_2q^{i+1},\alpha_2q,q;q)_i}\sum_{k=0}^i(-\alpha_1\alpha_2q^{(k+1)/2})^{-k}\frac{(\alpha_3 q^{n-i+1};q)_{i-k}(q^{-i},\alpha_1\alpha_2q^{i+1};q)_k}{(q;q)_k}f(k,0)
\end{equation}

\end{theorem}
\begin{proof}
Set

\[
S_i=q^{i(i+1)}\sum_{x_1=0}^i\frac{(\alpha_1q;q)_x(\alpha_2q;q)_{i-x_1}}{(q;q)_{x_1}(q;q)_{i-x_1}}(\alpha_1q)^{i-x_1}f(x_1,i-x_1)Q_i(x_1;\alpha_1,\alpha_2,i|q).
\]

\noindent
On one hand, $S_i$ may be seen as the scalar product of $f$ with $Q_i(x_1;\alpha_1,\alpha_2,i|q)$, along the line $x_1+x_2=i$. Hence using the orthogonality relations \eqref{qHahnortrel} for $Q_i(\cdot;\alpha_1,\alpha_2,i|q)$ and the formula $Q_{n-i}(0;\alpha_1\alpha_2q^{2i+1},\alpha_3,n-i|q)=q^{-(n-i)^2/2}(q;q)_{n-i}$ (see \eqref{Qn0}) we get:

\[
S_i=a_i\cdot\alpha_1^i\frac{(\alpha_1\alpha_2q^{i+1},\alpha_2 q,q;q)_i(q;q)_{n-i}}{(\alpha_1q;q)_i}q^{(n-i)i+(3i-n^2)/2}.
\] 

\noindent
On the other hand, applying \eqref{Duformula} and then \eqref{QnVandermonde} we get

\[
\begin{split}
S_i=&\sum_{k=0}^i(\alpha_1\alpha_2 q^{i+2})^{i-k}q^k(-1)^i\frac{(\alpha_1q;q)_k(\alpha_3q^{n-i+1};q)_{i-k}}{(q;q)_{i-k}}f(k,0)\sum_{x_1=0}^i(-1)^{x_1}\frac{Q_i(x_1;\alpha_1,\alpha_2,i|q)}{(q;q)_{x_1}(q;q)_{k-x_1}}q^{\frac{x_1(x_1-1)}{2}}\\
=&(-\alpha_1\alpha_2 q^{2+i/2})^i\sum_{k=0}^i(-\alpha_1\alpha_2q^{(k+1)/2})^{-k}\frac{(\alpha_3q^{n-i+1};q)_{i-k}(q^{-i},\alpha_1\alpha_2q^{i+1};q)_k}{(q;q)_k}f(k,0).
\end{split}
\]

\noindent
Comparing the two expressions for $S_i$ one gets immediately \eqref{ai}.

\end{proof}

Now we introduce the $q$-Racah polynomials giving them a particular normalization. More precisely, we set:

\[
r_n(x;\alpha,\beta,\delta,N|q)=q^{-n(N-n)}\frac{(\beta\delta q,q^{N-n+1};q)_n}{(\alpha\beta q^{n+1},q;q)_n}
\;_4\phi_3
\left[\begin{array}{c}
q^{-n},\delta q^{x-N},q^{-x},\alpha\beta q^{n+1}\\
\alpha q,\beta\delta q,q^{-N}\end{array};q,q
\right].
\]

\noindent
With respect to the standard definition (see \cite{GaRabook,Ismail,KoSw}) we have just added the factor $q^{-n(N-n)}\frac{(\beta\delta q,q^{N-n+1};q)_n}{(\alpha\beta q^{n+1},q;q)_n}$ (in the notation of \cite{Ismail,KoSw}, we also suppose that $\gamma q=q^{-N}$). If we denote by $\widetilde{r}_n$ the $q$-Racah polynomial in \cite{GaRa} then we have:

\begin{equation}\label{GaRaqRacah}
\widetilde{r}_n(x;\alpha,\beta,\delta,N|q)=\frac{(-1)^n(\alpha\beta q^{n+1},\alpha q,q;q)_n}{(q^{-N+n+1}\delta)^{n/2}}r_n(x;\alpha,\beta,\delta,N|q).
\end{equation}

We are ready to give the main result of this section, which is our version of Theorems 3.2 and 4.1 in \cite{Du6}.

\begin{theorem}\label{DunklqRacah}
For $j=0,1,2,\dotsc,n$ we have:

\begin{equation}\label{xirtheta}
\xi_{n-j,j}=\sum_{i=0}^n r_i(j;\alpha_2,\alpha_1,\alpha_2\alpha_3q^{n+1},n|q)\theta_{n-i,i}.
\end{equation}
\end{theorem}

\begin{proof}
In virtue of  4. in Theorem \ref{theomultqHahn}, it suffices to prove that

\[
\widetilde{
\xi}_{n-j,j}=\sum_{i=0}^n r_i(j;\alpha_2,\alpha_1,\alpha_2\alpha_3q^{n+1},n|q)\widetilde{\theta}_{n-i,i}.
\]

\noindent
Suppose that $\widetilde{\xi}_{n-j,j}=\sum_{i=0}^n a_{ij}\widetilde{\theta}_{n-i,i}$. Formulas \eqref{Qnn} and \eqref{Qn0} yields 

\[
\widetilde{\xi}_{n-j,j}(k,0)=(-\alpha_1)^kq^{k(k+1)/2-(n^2+j^2)/2}\frac{(q;q)_{n-j}(\alpha_2\alpha_3q^{n+j-k+2};q)_k(q^{n-j-k+1};q)_j}{(\alpha_1 q;q)_k}.
\]

\noindent
Therefore from \eqref{ai} we get

\[
\begin{split}
&a_{ij}=(-\alpha_2 q^{-n+(3i+1)/2})^i\frac{(\alpha_1 q;q)_i}{(q;q)_{n-i}(\alpha_1\alpha_2 q^{i+1},\alpha_2 q,q;q)_i}\sum_{k=0}^i
(-\alpha_2)^{-k}\frac{(\alpha_3 q^{n-i+1};q)_{i-k}}{(\alpha_1 q,q;q)_k}\\
\times &(q^{-i},\alpha_1\alpha_2 q^{i+1},\alpha_2\alpha_3 q^{n+j-k+2};q)_k(q^{n-j-k+1};q)_j(q;q)_{n-j}\\
&=(-\alpha_2 q^{-n+(3i+1)/2})^i\frac{(\alpha_1 q,\alpha_3 q^{n-i+1};q)_i(q^{n-i+1};q)_i}{(\alpha_1\alpha_2 q^{i+1},\alpha_2 q,q;q)_i}
\sum_{k=0}^i\frac{(q^{-i},\alpha_1\alpha_2 q^{i+1},\alpha_2^{-1}\alpha_3^{-1}q^{-n-j-1}, q^{-n+j};q)_k}{(q^{-n},\alpha_3^{-1}q^{-n},\alpha_1 q,q;q)_k}q^k\\
&=(-\alpha_2 q^{-n+(3i+1)/2})^i\frac{(\alpha_1 q,\alpha_3 q^{n-i+1};q)_i(q^{n-i+1};q)_i}{(\alpha_1\alpha_2 q^{i+1},\alpha_2 q,q;q)_i}
\;_4\phi_3
\left[\begin{array}{c}
q^{-i},\alpha_1\alpha_2 q^{i+1},\alpha_2^{-1}\alpha_3^{-1}q^{-n-j-1},q^{-n+j}\\
q^{-n},\alpha_3^{-1} q^{-n},\alpha_1 q\end{array};q,q
\right]\\
&=q^{-i(n-i)}\frac{(\alpha_1\alpha_2\alpha_3 q^{n+2},q^{n-i+1};q)_i}{(\alpha_1\alpha_2 q^{i+1},q;q)_i}
\;_4\phi_3
\left[\begin{array}{c}
q^{-i},\alpha_1\alpha_2 q^{i+1},\alpha_2\alpha_3q^{j+1},q^{-j}\\
q^{-n},\alpha_1\alpha_2\alpha_3 q^{n+2},\alpha_2 q\end{array};q,q
\right]\\
&\equiv  r_i(j;\alpha_2,\alpha_1,\alpha_2\alpha_3q^{n+1},n|q),
\end{split}
\]

\noindent
where in the second equality we have used several times the elementary identities $(a^{-1};q)_h=a^{-h}(-1)^h$ $q^{h(h-1)/2}(aq^{-h+1};q)_h$ and $(aq^{m-h+1};q)_{h-r}(aq^{m-r+1};q)_r=(aq^{m-h+1};q)_h$ and in the fourth equality we have used the Sears $\;_4\phi_3$ transformation formula (2.10.4) in \cite{GaRabook}, with $n=i$, $a=\alpha_1\alpha_2q^{i+1}$, $b=\alpha_2^{-1}\alpha_3^{-1}q^{-n-j-1}$, $c=q^{j-n}$, $d=q^{-n}$, $e=\alpha_3^{-1}q^{-n}$ and $f=\alpha_2 q$.
\end{proof}

%%%%%%%%%%%%%%%%%%%%%%%%%%%%%%%%%%%%%%%%%%%%%%%%%%%%%%%%%%%%%%%%%%%%%%%%%%%%%%%%%%
%%%%%%%%%%%%%%%%%%%%%%%%%%%%%%%%%%%%%%%%%%%%%%%%%%%%%%%%%%%%%%%%%%%%%%%%%%%%%%%%%%%%%
\section{Connection coefficients for the transplantation of an edge}\label{sectiontransplantation}
%%%%%%%%%%%%%%%%%%%%%%%%%%%%%%%%%%%%%%%%%%%%%%%%%%%%%%%%%%%%%%%%%%%%%%%%%%%%%%%%%%%%%

Now we show that the $q$-Racah polynomials are also the connection coefficients for a basic operation called the {\em transplantation} of an edge. Let $h,r,s$ be three positive integer with $1\leq s<r<h$. Suppose that $\mathcal{T}'$, $\mathcal{T}''$ and $\mathcal{T}'''$ rooted binary trees with respectively $s$, $r-s$, $h-r$ leaves. In the case $s=1$ then we suppose that $\mathcal{T}'$ is a single vertex; similarly if $r-s=1$ or $h-r=1$. Then we can construct two different rooted trees with $h$ leaves as in the figures below:

\begin{picture}(400,90)
%%%%%%%%%%%%%%%%%%%%%%%%%%%%%%%%%%%
%%%%%%%%%%%%%%% first tree
%%%%%%%%%%%%%%%%%%%%%%%%%%%%%%%
\thicklines
%%%%%%%%%%%%%%% zero level
\put(100,70){\circle*{4}}
\put(100,74){$V$}
\put(104,40){$a$}

%%%%%%%%%%%% first level 
\put(80,50){\circle*{4}}
\put(120,50){\circle*{4}}

\thicklines
\put(100,70){\line(-1,-1){20}}
\put(100,70){\line(1,-1){20}}
\drawline(120,50)(140,30)

%%%%%%%%%%%%%%%  second level
\put(100,30){\circle*{4}}
\put(140,30){\circle*{4}}

%%%%%%%%%% labelings
\put(74,36){$\mathcal{T}'$}
\put(93,16){$\mathcal{T}''$}
\put(132,16){$\mathcal{T}'''$}

\put(120,50){\line(-1,-1){20}}

\thinlines
\put(80,50){\line(-1,-1){10}}
\put(80,50){\line(1,-1){10}}

\put(100,30){\line(-1,-1){10}}
\put(100,30){\line(1,-1){10}}

\put(140,30){\line(-1,-1){10}}
\put(140,30){\line(1,-1){10}}

\put(78,0){{\bf Figure 4a}}

\put(200,55){\large $\stackrel{\tau}{\longrightarrow}$}

%%%%%%%%%%%%%%% second tree

%%%%%%%%%%%%%%% zero level
\put(320,70){\circle*{4}}
\put(320,74){$V$}
\put(313,40){$a$}

%%%%%%%%%%%% first level 
\put(340,50){\circle*{4}}
\put(300,50){\circle*{4}}

\thicklines
\put(320,70){\line(-1,-1){20}}
\put(320,70){\line(1,-1){20}}
\drawline(300,50)(280,30)

%%%%%%%%%% labelings
%\put(87,105){$Z$}
%\put(45,105){$W$}

%%%%%%%%%%%%%%%  second level
\put(280,30){\circle*{4}}
\put(320,30){\circle*{4}}
\put(300,50){\line(1,-1){20}}

\thinlines
\put(340,50){\line(-1,-1){10}}
\put(340,50){\line(1,-1){10}}

\put(280,30){\line(-1,-1){10}}
\put(280,30){\line(1,-1){10}}

\put(320,30){\line(-1,-1){10}}
\put(320,30){\line(1,-1){10}}

%%%%%%%%%%%% labelings
\put(313,15){$\mathcal{T}''$}
\put(275,15){$\mathcal{T}'$}
\put(332,35){$\mathcal{T}'''$}

\put(275,0){{\bf Figure 4b}}
\end{picture}
\quad\\
\quad\\

\noindent
The passage from the tree in Figure 4a to the tree in Figure 4b is called the {\em transplantation} $\tau$ of the edge $a$ from {\em right to left}. With respect to the standard theory exposed in Section 10.5 of \cite{K-V} and Section 6.3 of \cite{NSU} (see also \cite{VdJ}), we do not consider the operation of transpositions: our equations are not invariant under permutations of coordinate. Moreover, we will essentially consider only transplantations from right to left, as above; see Remark \ref{righttoleft} for the reason. Clearly we may consider transplantations of edges at every level: the root $V$ may be replaced by any vertex $U$ of the tree and in this case we will call of a transplantation of an edge immediately below $U$. \\

Now fix a positive integer $N$ and let $\alpha_1,\alpha_2,\dotsc,\alpha_h$ be real numbers satisfying the conditions $0<\alpha_i<q^{-1}$, $i=1,2,\dotsc,h$, or the conditions $\alpha_i>q^{-N}$, $i=1,2,\dotsc,h$. Fix also coefficients labellings $\text{c}'$, $\text{c}''$ and $\text{c}'''$ respectively of $\mathcal{T}'$, $\mathcal{T}''$ and $\mathcal{T}'''$. Set ${\bf x}'=(x_1,x_2,\dotsc,x_s)$, ${\bf x}''=(x_{s+1},x_{s+2},\dotsc,x_r)$ and ${\bf x}'''=(x_{r+1},x_{r+2},\dotsc,x_h)$. We can consider the multidimensional $q$-Hahn polynomials associated to those labeled trees: 
$Q_{\text{c}'}({\bf x}';\alpha_1,\dotsc,\alpha_s,x'|q)$, $Q_{\text{c}''}({\bf x}'';\alpha_{s+1},\dotsc,\alpha_r,x''|q)$ and $Q_{\text{c}'''}({\bf x}''';\alpha_{r+1},\dotsc,\alpha_h,x'''|q)$, that we denote simply by $Q_{\text{c}'}$, $Q_{\text{c}''}$ and $Q_{\text{c}'''}$.
Fix also a positive integer $n$ and nonnegative integers $u,v,i,l,j$ such that: $u+j\leq n$, $i+l\leq u$, $v+i\leq n$, $l+j\leq v$. Then taking $\text{c}'\in\text{CL}(\mathcal{T}',i)$, $\text{c}''\in\text{CL}(\mathcal{T}'',l)$ and $\text{c}'''\in\text{CL}(\mathcal{T}''',j)$ we can label the threes in Figures 4a and 4b as in the figures below:

\begin{picture}(400,90)
%%%%%%%%%%%%%%% first tree

\thicklines
%%%%%%%%%%%%%%% zero level
\put(100,70){\circle*{4}}
\put(85,74){$n-v-i$}

%%%%%%%%%%%% first level 
\put(80,50){\circle*{4}}
\put(120,50){\circle*{4}}
\put(124,52){$v-l-j$}

\thicklines
\put(100,70){\line(-1,-1){20}}
\put(100,70){\line(1,-1){20}}
\drawline(120,50)(140,30)

%%%%%%%%%%%%%%%  second level
\put(100,30){\circle*{4}}
\put(140,30){\circle*{4}}

%%%%%%%%%% labelings
\put(76,38){$\text{c}'$}
\put(95,17){$\text{c}''$}
\put(134,16){$\text{c}'''$}

\put(120,50){\line(-1,-1){20}}

\thinlines
\put(80,50){\line(-1,-1){10}}
\put(80,50){\line(1,-1){10}}

\put(100,30){\line(-1,-1){10}}
\put(100,30){\line(1,-1){10}}

\put(140,30){\line(-1,-1){10}}
\put(140,30){\line(1,-1){10}}

\put(78,0){{\bf Figure 5a}}

%%%%%%%%%%%%%%%%%%%%%%%%%%%%%%%%%%%
%%%%%%%%%%%%%%% second tree
%%%%%%%%%%%%%%%%%%%%%%%%%%%%%%%
%%%%%%%%%%%%%%% zero level

\put(320,70){\circle*{4}}
%%%%%% labelings
\put(305,74){$n-u-j$}

%%%%%%%%%%%% first level 
\put(340,50){\circle*{4}}
\put(300,50){\circle*{4}}

\thicklines
\put(320,70){\line(-1,-1){20}}
\put(320,70){\line(1,-1){20}}
\drawline(300,50)(280,30)

\put(300,50){\line(1,-1){20}}

\thinlines
\put(340,50){\line(-1,-1){10}}
\put(340,50){\line(1,-1){10}}

\put(280,30){\line(-1,-1){10}}
\put(280,30){\line(1,-1){10}}

\put(320,30){\line(-1,-1){10}}
\put(320,30){\line(1,-1){10}}

%%%%%%%%%% labelings
\put(260,52){$u-i-l$}

%%%%%%%%%%%%%%%  second level
\put(280,30){\circle*{4}}
\put(320,30){\circle*{4}}

%%%%%%%%%%%% labelings
\put(314,17){$\text{c}''$}
\put(276,18){$\text{c}'$}
\put(333,35){$\text{c}'''$}
\put(275,0){{\bf Figure 5b}}

\end{picture}
\quad\\

\noindent
Denote by $\Xi_v(\bf{x})$ and $\Theta_u({\bf x})$ the multidimensional $q$-Hahn polynomials associated respectively to the tree in Figure 4a and 4b, with the parameters $\alpha_1,\alpha_2,\dotsc,\alpha_h$ and coefficients as in Figure 5a and 5b (we allow $v$ and $u$ to vary). From the definition \eqref{multqHahn} we have:

\[
\begin{split}
\Xi_v({\bf x})= &q^{-vX_s}Q_{n-v-i}(X_s-i;A_s q^{s+2i-1}, A_s^{-1}A_h q^{h-s+2v-1},X_h-v-i|q)  q^{-j(X_r-X_s)}\\
&\times Q_{v-l-j}(X_r-X_s-l;A_s^{-1}A_r q^{r-s+2l-1}, A_r^{-1}A_h q^{h-r+2j-1}, X_h-X_s-l-j|q)Q_{\text{c}'}Q_{\text{c}''}Q_{\text{c}'''}\\
\equiv & q^{-i(v-l-j)-lX_s-jX_r}\xi_{n-v-i,v-l-j}(X_s-i,X_r-X_s-l,X_h-X_r-j;A_s q^{s+2i-1},\\
&A_s^{-1}A_r q^{r-s+2l-1},
A_r^{-1}A_h q^{h-r+2j-1})Q_{\text{c}'}Q_{\text{c}''}Q_{\text{c}'''},
\end{split}
\] 

\noindent
and similarly

\[
\begin{split}
\Theta_u({\bf x})=&q^{-jX_r}Q_{n-u-j}(X_r-u;A_r q^{r+2u-1},A_r^{-1}A_h q^{h-r+2j-1},X_h-u-j|q)\\
&\times q^{-lX_s}Q_{u-i-l}(X_s-i;A_s q^{s+2i-1},A_s^{-1}A_r q^{r-s+2l-1},X_r-i-l|q)Q_{\text{c}'}Q_{\text{c}''}Q_{\text{c}'''}\\
\equiv & q^{-lX_s-jX_r}\theta_{n-u-j,u-i-l}(X_s-i,X_r-X_s-l,X_h-X_r-j;A_s q^{s+2i-1},A_s^{-1}A_r q^{r-s+2l-1},\\
&\quad A_r^{-1}A_h q^{h-r+2j-1})Q_{\text{c}'}Q_{\text{c}''}Q_{\text{c}'''}
\end{split}
\] 

\noindent
where $\theta_{n-u-j,u-i-l}$ and $\xi_{n-v-i,v-l-j}$ are as in Section \ref{bidqHahn}. Therefore, from \eqref{xirtheta} it follows immediately the following {\em connection formula} for the transplantation of an edge:

\begin{equation}\label{XirTheta}
\begin{split}
\Xi_v=q^{-i(v-l-j)}\sum_{u=i+l}^{n-j}
r_{u-i-l}(&v-l-j;A_s^{-1}A_r q^{r-s+2l-1},A_s q^{s+2i-1},\\
&A_s^{-1}A_h q^{n+l+j-i+h-s-1},n-i-l-j|q)\Theta_u,
\end{split}
\end{equation}

\begin{remark}\label{righttoleft}
{\rm

We have arranged the normalizations of the $q$-Hahn and $q$-Racah polynomials in order to get the formulas \eqref{xirtheta} and \eqref{XirTheta} as simple as possible. But this leads to more complicated formulas for  transplantations from {\em left to right}; compare with Theorem 3.2 and Corollary 3.3 in \cite{Du6}. Therefore, for simplicity, we will consider only transplantations from right to left. In any case, this way we can construct explicitly a wide class of multidimensional $q$-Racah polynomials.

}
\end{remark}

%%%%%%%%%%%%%%%%%%%%%%%%%%%%%%%%%%%%%%%%%%%%%%%%%%%%%%%%%%%%%%%%%%%%%%%%%%%%%%%%%%
%%%%%%%%%%%%%%%%%%%%%%%%%%%%%%%%%%%%%%%%%%%%%%%%%%%%%%%%%%%%%%%%%%%%%%%%%%%%%%%%%%%%%
\section{Multidimensional $q$-Racah polynomials }\label{sectionmultqRacah}
%%%%%%%%%%%%%%%%%%%%%%%%%%%%%%%%%%%%%%%%%%%%%%%%%%%%%%%%%%%%%%%%%%%%%%%%%%%%%%%%%%%%%

We formulate and prove a particular case of a well known result (\cite{K-V}, \cite{NSU}).

\begin{proposition}\label{sequencetranspl}
Let $\mathcal{T}$ and $\mathcal{S}$ be two rooted binary trees with $h$ leaves.

\begin{enumerate}

\item If $\mathcal{S}$ is as in Example \ref{bidqHahn2} then there exists a sequence of transplantations from right to left that transforms $\mathcal{T}$ in $\mathcal{S}$. 

\item If $\mathcal{T}$ is as in Example \ref{bidqHahn1} then there exists a sequence of transplantations from right to left that transforms $\mathcal{T}$ in $\mathcal{S}$. 

\end{enumerate}
\end{proposition}
\begin{proof}
For instance, in the first case we can perform repeatedly the highest possible transplantation from right to left. This procedure reduces an arbitrary $\mathcal{T}$ to the tree in Example \ref{bidqHahn2}.  
\end{proof}

Now we are in position to define our multidimensional $q$-Racah polynomials. Let $h$, $N$, $q$ and $\alpha_1,\alpha_2,\dotsc,\alpha_h$ be as in the previous sections and let $\mathcal{T}$ and $\mathcal{S}$ be two rooted binary trees with $h$ leaves. 
From 2. in Theorem \ref{theomultqHahn} (see also Theorem \ref{DunklqRacah} and its proof) we deduce that to determine the connection coefficients between the basis for $V_{h,N}$ associated to the trees $\mathcal{T}$ and $\mathcal{S}$ it suffices to analyze the case $n=N$. For $\text{c}\in\text{CL}(\mathcal{T},n)$ and $\text{d}\in\text{CL}(\mathcal{S},n)$ we set 

\[
Q_\text{c}=Q_\text{c}(\cdot;\alpha_1,\alpha_2,\dotsc,\alpha_h|q)\qquad\qquad\text{and}\qquad\qquad Q_\text{d}=Q_\text{d}(\cdot;\alpha_1,\alpha_2,\dotsc,\alpha_h|q).
\]

\noindent
We define the {\em multidimensional} $q$-{\em Racah polynomials associated to the trees} $\mathcal{T}$ and $\mathcal{S}$ as the connection coefficients $\{r_\text{d}(\text{c}):\text{d}\in\text{CL}(\mathcal{S},n),\text{c}\in\text{CL}(\mathcal{T},n)\}$ between the bases $\{Q_\text{c}:\text{c}\in\text{CL}(\mathcal{T},n)\}$ and $\{Q_\text{d}:\text{d}\in\text{CL}(\mathcal{S},n)\}$ of $\text{Ker}\mathcal{L}_n$. In formul\ae:

\[
Q_\text{c}=\sum_{\text{d}\in\text{CL}(\mathcal{S},n)} r_\text{d}(\text{c})Q_\text{d}
\]

\noindent
The coefficients $r_\text{d}(\text{c})$ satisfy the following orthogonality relations:

\begin{equation}\label{orthrelmultqRacah}
\sum_{\text{c}\in\text{CL}(\mathcal{T},n)} r_{\text{d}_1}(\text{c})r_{\text{d}_2}(\text{c})\frac{1}{\lVert Q_\text{c}\rVert^2}=\frac{\delta_{\text{d}_1,\text{d}_2}}{\lVert Q_{\text{d}_1}\rVert^2}.
\end{equation} 

Now suppose that $\tau_1,\tau_2,\dotsc,\tau_p$ a sequence of transplantations from right to left that transfoms $\mathcal{T}$ in $\mathcal{S}$.
As an immediate consequence of Proposition \ref{sequencetranspl}, we get that for arbitrary $\mathcal{T}$ and $\mathcal{S}$ there exists a sequence of transplantations that transforms $\mathcal{T}$ in $\mathcal{S}$, but this sequence might contain also transplantations from left to right. Therefore, in our assumptions, the couple $\mathcal{T}$ and $\mathcal{S}$ is not arbitrary (see also Remark \ref{righttoleft}). If we fix also $\text{c}\in\text{CL}(\mathcal{T},n)$ and $\text{d}\in\text{CL}(\mathcal{S},n)$ then 
for each $\tau_k$, $k=1,2,\dotsc,p$, there  exist a set of {\em transplantation coefficients} $s_k, r_k, n_k, i_k, l_k, j_k, u_k, v_k$ as in Figures 5a, 5a and equation \eqref{XirTheta} that lead from $\mathcal{T}$ with labeling $\text{c}$ to $\mathcal{S}$ with labeling $\text{d}$. Therefore we get the following explicit formula (or explicit algorithm) to compute $r_\text{d}(\text{c})$:

\begin{multline}\label{multqRacahexpression}
r_\text{d}(\text{c})=\prod_{k=1}^p q^{-i_k(v_k-l_k-j_k)}
r_{u_k-i_k-l_k}(v_k-l_k-j_k;A_{s_k}^{-1}A_{r_k} q^{r_k-s_k+2l_k-1},A_{s_k} q^{s_k+2i_k-1},\\
A_{s_k}^{-1}A_h q^{n+l_k+j_k-i_k+h-s_k-1},n-i_k-l_k-j_k|q)
\end{multline}

\begin{example}{\rm

Consider the following sequence of transplantations, that starts with the tree in Example \ref{bidqHahn1} for $h=5$:

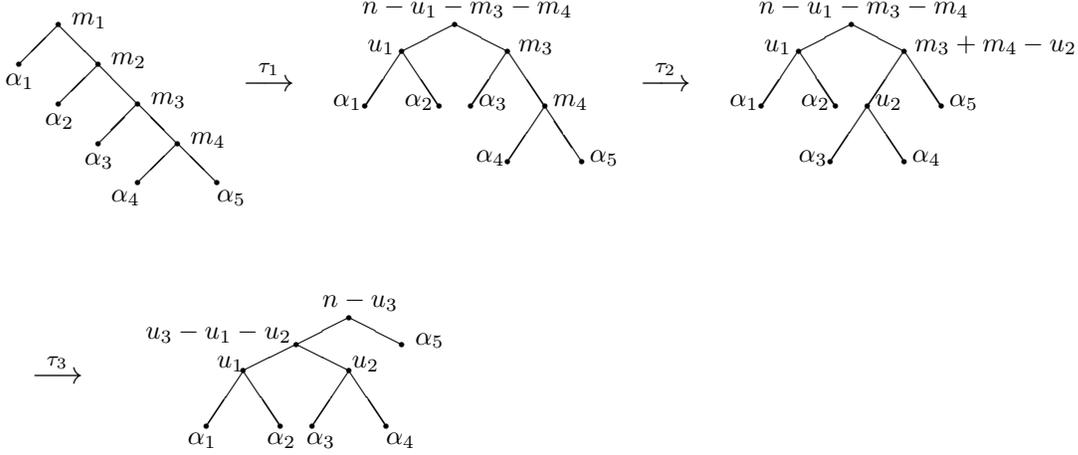
\begin{picture}(400,100)

%%%%%%%%%%%%%%%%%%%%%%%%%%%%
%%%%%%%%%%%%%%% first tree
%%%%%%%%%%%%%%%%%%%%%%%%%%%%
\thinlines
%%%%%%%%%%%%%%% zero level
\put(15,80){\circle*{2}}

%%%%%%%%%%%% first level 

\put(15,80){\line(-1,-1){15}}
\put(15,80){\line(1,-1){30}}
\put(0,65){\circle*{2}}
\put(30,65){\circle*{2}}

%%%%%% labelings
\put(20,80){$m_1$}

%%%%%%%%%% labelings
\put(-5,57){$\alpha_1$}
\put(35,65){$m_2$}
\put(85,55){\large $\stackrel{\tau_1}{\longrightarrow}$}

%%%%%%%%%%%%%%%  second level
\put(30,65){\line(-1,-1){15}}

\put(15,50){\circle*{2}}
\put(45,50){\circle*{2}}

%%%%%%%%%%%% labelings
\put(10,42){$\alpha_2$}
\put(50,50){$m_3$}

%%%%%%%%%%%%%%%%%%% third level
\drawline(45,50)(60,35)
\put(45,50){\line(-1,-1){15}}
\put(60,35){\circle*{2}}
\put(30,35){\circle*{2}}

%%%%%%%%%%%%%%%% labelings
\put(25,27){$\alpha_3$}
\put(65,35){$m_4$}

%%%%%%%%%%%%%%%%%%%%%%%% last level
\put(60,35){\line(-1,-1){15}}
\put(60,35){\line(1,-1){15}}
\put(75,20){\circle*{2}}
\put(45,20){\circle*{2}}

\put(35,13){$\alpha_4$}
\put(75,13){$\alpha_5$}

%%%%%%%%%%%%%%%%%%%%%%%%%%%%
%%%%%%%%%%%%%%% second tree
%%%%%%%%%%%%%%%%%%%%%%%%%%%%%

%%%%%%%%%%%%%%% zero level
\put(165,80){\circle*{2}}
\put(130,84){$n-u_1-m_3-m_4$}

%%%%%%%%%%%% first level 
\put(165,80){\line(-2,-1){20}}
\put(165,80){\line(2,-1){20}}
\put(145,70){\circle*{2}}
\put(185,70){\circle*{2}}

%%%%%% labelings
\put(132,70){$u_1$}
\put(189,70){$m_3$}
\put(235,55){\large $\stackrel{\tau_2}{\longrightarrow}$}

%%%%%%%%%%%% second level
\put(145,70){\line(-2,-3){14}}
\put(145,70){\line(2,-3){14}}
\put(131,49){\circle*{2}}
\put(159,49){\circle*{2}}

\put(185,70){\line(-2,-3){14}}
\put(185,70){\line(2,-3){14}}
\put(171,49){\circle*{2}}
\put(199,49){\circle*{2}}

%%%%%%%%%%%% labelings
\put(119,49){$\alpha_1$}
\put(146,49){$\alpha_2$}
\put(174,49){$\alpha_3$}
\put(202,49){$m_4$}

%%%%%%%%%%%% third level
\put(199,49){\line(-2,-3){14}}
\drawline(199,49)(213,28)
\put(185,28){\circle*{2}}
\put(213,28){\circle*{2}}

%%%%%%%%%%%% labelings
\put(173,28){$\alpha_4$}
\put(216,28){$\alpha_5$}

%%%%%%%%%%%%%%%%%%%%%%%%%%%%
%%%%%%%%%%%%%%% third tree
%%%%%%%%%%%%%%%%%%%%%%%%%%%%%

%%%%%%%%%%%%%%% zero level
\put(315,80){\circle*{2}}
\put(280,84){$n-u_1-m_3-m_4$}

%%%%%%%%%%%% first level 
\put(315,80){\line(-2,-1){20}}
\put(315,80){\line(2,-1){20}}
\put(295,70){\circle*{2}}
\put(335,70){\circle*{2}}

%%%%%% labelings
\put(282,70){$u_1$}
\put(339,70){$m_3+m_4-u_2$}

%%%%%%%%%%%% second level
\put(295,70){\line(-2,-3){14}}
\put(295,70){\line(2,-3){14}}
\put(281,49){\circle*{2}}
\put(309,49){\circle*{2}}

\put(335,70){\line(-2,-3){14}}
\put(335,70){\line(2,-3){14}}
\put(321,49){\circle*{2}}
\put(349,49){\circle*{2}}

%%%%%%%%%%%% labelings
\put(269,49){$\alpha_1$}
\put(296,49){$\alpha_2$}
\put(324,49){$u_2$}
\put(352,49){$\alpha_5$}

%%%%%%%%%%%% third level
\put(321,49){\line(2,-3){14}}
\put(321,49){\line(-2,-3){14}}
\put(307,28){\circle*{2}}
\put(335,28){\circle*{2}}

%%%%%%%%%%%% labelings
\put(295,28){$\alpha_3$}
\put(338,28){$\alpha_4$}

\end{picture}

\begin{picture}(400,90)

\put(5,35){\large $\stackrel{\tau_3}{\longrightarrow}$}

%%%%%%%%%%%%%%% zero level
\put(125,60){\circle*{2}}
\put(115,64){$n-u_3$}

%%%%%%%%%%%% first level 
\put(125,60){\line(-2,-1){40}}
\put(125,60){\line(2,-1){20}}
\put(145,50){\circle*{2}}
\put(105,50){\circle*{2}}

%%%%%% labelings
\put(48,52){$u_3-u_1-u_2$}
\put(150,50){$\alpha_5$}

%%%%%%%%%%%% second level
\put(105,50){\line(2,-1){20}}
\put(85,40){\circle*{2}}
\put(125,40){\circle*{2}}

%%%%%%%%%%%% labelings
\put(75,41){$u_1$}
\put(126,41){$u_2$}

%%%%%%%%%%%% third level
\put(125,40){\line(2,-3){14}}
\put(125,40){\line(-2,-3){14}}
\put(139,19){\circle*{2}}
\put(111,19){\circle*{2}}

\put(85,40){\line(2,-3){14}}
\put(85,40){\line(-2,-3){14}}
\put(99,19){\circle*{2}}
\put(71,19){\circle*{2}}

%%%%%%%%%%%% labelings
\put(64,12){$\alpha_1$}
\put(94,12){$\alpha_2$}
\put(109,12){$\alpha_3$}
\put(139,12){$\alpha_4$}

\end{picture}

\noindent
The associated multidimensional $q$-Racah polynomials is the following:

\begin{multline}\label{3dimqRacah}
r_{u_1,u_2,u_3}(m_1,m_2,m_3,m_4|q)=\\
r_{u_1}(m_2;\alpha_2,\alpha_1,\alpha_2\alpha_3\alpha_4\alpha_5 q^{n+m_3+m_4+3},m_1+m_2|q)\cdot r_{u_2}(m_4;\alpha_4,\alpha_3,\alpha_4\alpha_5 q^{m_3+m_4+1},m_3+m_4|q)\\
\times q^{-u_1(m_3+m_4-u_2)}r_{u_3-u_1-u_2}(m_3+m_4-u_2;\alpha_3\alpha_4q^{2u_2+1},\alpha_1\alpha_2 q^{2u_1+1},\alpha_3\alpha_4\alpha_5 q^{n+u_2-u_1+2},n-u_1-u_2|q).
\end{multline}

\noindent
We note that taking the limit $\alpha_5\rightarrow 0$ and using the known relation between $q$-Racah and $q$-Hahn polynomials \cite{GaRabook,GaRa,KoSw}, the three dimensional $q$-Racah polynomial \eqref{3dimqRacah} becomes a multiple of the three dimensional $q$-Hahn polynomial in Example \ref{3dimqHahn}. If we take \eqref{3dimqRacah} for $0\leq u_3\leq n$, $u_1,u_2\geq 0$ and $u_1+u_2\leq u_3$, we get a complete family of orthogonal polynomials on the set $\{(m_1,m_2,m_3,m_4):m_1+m_2+m_3+m_4=n,m_i\geq 0,i=1,2,3,4\}$. From \eqref{orthrelmultqRacah} we know that the weight in the orthogonality relations is just the reciprocal of \eqref{normxi} for $N=n$ and $h=5$, while the square of the norm is equal to the reciprocal of the square of the norm of the $q$-Hahn polynomials associated to the last tree (with $n=N$), and therefore from \eqref{normQc} is given by:

\[
\begin{split}
&\frac{(A_4q^{2u_3+4};q)_{n-u_3}(\alpha_1\alpha_2 q^{2u_1+2};q)_{u_3-u_1-u_2}(\alpha_1q;q)_{u_1}}{(q,A_5q^{n+u_3+4},\alpha_5q;q)_{n-u_3}(q,A_4q^{u_1+u_2+u_3+3},\alpha_3\alpha_4q^{2u_2+2})_{u_3-u_1-u_2}(q,\alpha_1\alpha_2q^{u_1+1},\alpha_2q;q)_{u_1}}\\
&\times \frac{(\alpha_3q;q)_{u_2}\alpha_1^{-n}\alpha_2^{-n+u_1}\alpha_3^{-n+u_3-u_2}\alpha_4^{-n+u_3}}{(q,\alpha_3\alpha_4q^{u_2+1},\alpha_4q;q)_{u_2}}q^{-(2u_3+3)(n-u_3)-(2u_1+1)(u_3-u_1)+2u_1u_2-u_2+(n^2-3n)/2}.\\
\end{split}
\]

}
\end{example}

%%%%%%%%%%%%%%%%%%%%%%%%%%%%%%%%%%%%%%%%%%%%%%%%%%%%%%%%%%%%%%%%%%%%%%%%%%%%%%%%%%
%%%%%%%%%%%%%%%%%%%%%%%%%%%%%%%%%%%%%%%%%%%%%%%%%%%%%%%%%%%%%%%%%%%%%%%%%%%%%%%%%%%%%
\section{The polynomials of Gasper and Rahman}\label{sectionGaRa}
%%%%%%%%%%%%%%%%%%%%%%%%%%%%%%%%%%%%%%%%%%%%%%%%%%%%%%%%%%%%%%%%%%%%%%%%%%%%%%%%%%%%%

In this section we compute the connection coefficients between the $q$-Hahn polynomials in Examples \ref{bidqHahn1} and \ref{bidqHahn2} and show that they coincide with the multidimensional $q$-Racah polynomials in \cite{GaRa}. A sequence of transplantations that leads from the labeled tree in Example \ref{bidqHahn1} to the labeled tree in Example \eqref{bidqHahn2} is the following:

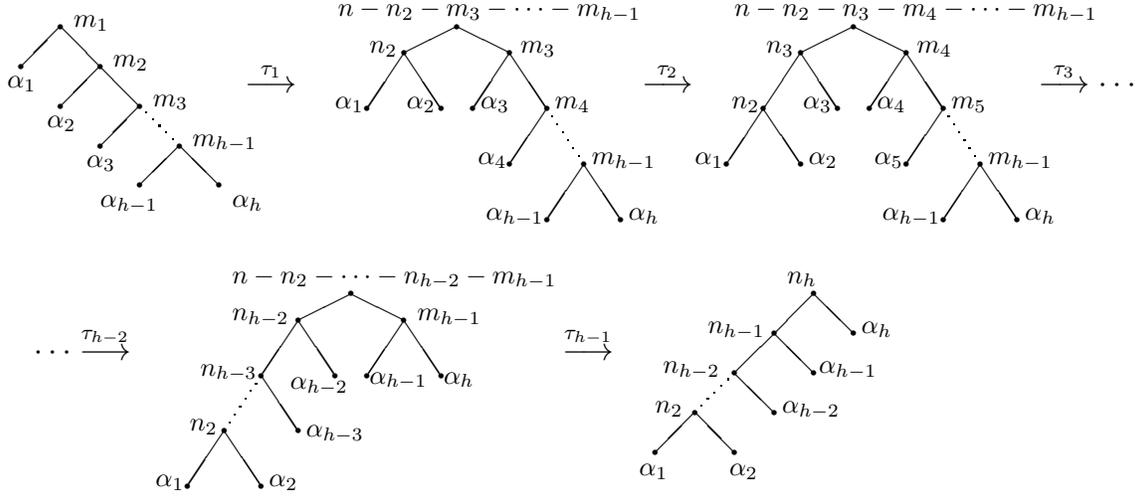
\begin{picture}(400,100)

%%%%%%%%%%%%%%%%%%%%%%%%%%%%
%%%%%%%%%%%%%%% first tree
%%%%%%%%%%%%%%%%%%%%%%%%%%%%
\thinlines
%%%%%%%%%%%%%%% zero level
\put(15,80){\circle*{2}}

%%%%%%%%%%%% first level 

\put(15,80){\line(-1,-1){15}}
\put(15,80){\line(1,-1){30}}
\put(0,65){\circle*{2}}
\put(30,65){\circle*{2}}

%%%%%% labelings
\put(20,80){$m_1$}

%%%%%%%%%% labelings
\put(-5,57){$\alpha_1$}
\put(35,65){$m_2$}
\put(85,55){\large $\stackrel{\tau_1}{\longrightarrow}$}

%%%%%%%%%%%%%%%  second level
\put(30,65){\line(-1,-1){15}}

\put(15,50){\circle*{2}}
\put(45,50){\circle*{2}}

%%%%%%%%%%%% labelings
\put(10,42){$\alpha_2$}
\put(50,50){$m_3$}

%%%%%%%%%%%%%%%%%%% third level
\dottedline{4}(45,50)(60,35)
\put(45,50){\line(-1,-1){15}}
\put(60,35){\circle*{2}}
\put(30,35){\circle*{2}}

%%%%%%%%%%%%%%%% labelings
\put(25,27){$\alpha_3$}
\put(65,35){$m_{h-1}$}

%%%%%%%%%%%%%%%%%%%%%%%% last level
\put(60,35){\line(-1,-1){15}}
\put(60,35){\line(1,-1){15}}
\put(75,20){\circle*{2}}
\put(45,20){\circle*{2}}

\put(30,12){$\alpha_{h-1}$}
\put(80,12){$\alpha_h$}

%%%%%%%%%%%%%%%%%%%%%%%%%%%%
%%%%%%%%%%%%%%% second tree
%%%%%%%%%%%%%%%%%%%%%%%%%%%%%

%%%%%%%%%%%%%%% zero level
\put(165,80){\circle*{2}}
\put(120,84){$n-n_2-m_3-\dotsb-m_{h-1}$}

%%%%%%%%%%%% first level 
\put(165,80){\line(-2,-1){20}}
\put(165,80){\line(2,-1){20}}
\put(145,70){\circle*{2}}
\put(185,70){\circle*{2}}

%%%%%% labelings
\put(132,70){$n_2$}
\put(189,70){$m_3$}
\put(235,55){\large $\stackrel{\tau_2}{\longrightarrow}$}

%%%%%%%%%%%% second level
\put(145,70){\line(-2,-3){14}}
\put(145,70){\line(2,-3){14}}
\put(131,49){\circle*{2}}
\put(159,49){\circle*{2}}

\put(185,70){\line(-2,-3){14}}
\put(185,70){\line(2,-3){14}}
\put(171,49){\circle*{2}}
\put(199,49){\circle*{2}}

%%%%%%%%%%%% labelings
\put(119,49){$\alpha_1$}
\put(146,49){$\alpha_2$}
\put(174,49){$\alpha_3$}
\put(202,49){$m_4$}

%%%%%%%%%%%% third level
\put(199,49){\line(-2,-3){14}}
\dottedline{4}(199,49)(213,28)
\put(185,28){\circle*{2}}
\put(213,28){\circle*{2}}

%%%%%%%%%%%% labelings
\put(173,28){$\alpha_4$}
\put(216,28){$m_{h-1}$}

%%%%%%%%%%%% last level
\put(213,28){\line(-2,-3){14}}
\put(213,28){\line(2,-3){14}}
\put(227,7){\circle*{2}}
\put(199,7){\circle*{2}}

%%%%%%%%%%%% labelings
\put(176,7){$\alpha_{h-1}$}
\put(230,7){$\alpha_h$}

%%%%%%%%%%%%%%%%%%%%%%%%%%%%
%%%%%%%%%%%%%%% third tree
%%%%%%%%%%%%%%%%%%%%%%%%%%%%%

%%%%%%%%%%%%%%% zero level
\put(315,80){\circle*{2}}
\put(270,84){$n-n_2-n_3-m_4-\dotsb-m_{h-1}$}

%%%%%%%%%%%% first level 
\put(315,80){\line(-2,-1){20}}
\put(315,80){\line(2,-1){20}}
\put(295,70){\circle*{2}}
\put(335,70){\circle*{2}}

%%%%%% labelings
\put(282,70){$n_3$}
\put(339,70){$m_4$}
\put(385,55){\large $\stackrel{\tau_3}{\longrightarrow}\dotsb$}

%%%%%%%%%%%% second level
\put(295,70){\line(-2,-3){28}}
\put(295,70){\line(2,-3){14}}
\put(281,49){\circle*{2}}
\put(309,49){\circle*{2}}

\put(335,70){\line(-2,-3){14}}
\put(335,70){\line(2,-3){14}}
\put(321,49){\circle*{2}}
\put(349,49){\circle*{2}}

%%%%%%%%%%%% labelings
\put(269,49){$n_2$}
\put(296,49){$\alpha_3$}
\put(324,49){$\alpha_4$}
\put(352,49){$m_5$}

%%%%%%%%%%%% third level
\put(349,49){\line(-2,-3){14}}
\dottedline{4}(349,49)(363,28)
\put(335,28){\circle*{2}}
\put(363,28){\circle*{2}}
\put(281,49){\line(2,-3){14}}
\put(267,28){\circle*{2}}
\put(295,28){\circle*{2}}

%%%%%%%%%%%% labelings
\put(323,28){$\alpha_5$}
\put(366,28){$m_{h-1}$}
\put(255,28){$\alpha_1$}
\put(298,28){$\alpha_2$}

%%%%%%%%%%%% last level
\put(363,28){\line(-2,-3){14}}
\put(363,28){\line(2,-3){14}}
\put(377,7){\circle*{2}}
\put(349,7){\circle*{2}}

%%%%%%%%%%%% labelings
\put(326,7){$\alpha_{h-1}$}
\put(380,7){$\alpha_h$}

\end{picture}

\begin{picture}(400,100)

%%%%%%%%%%%%%%%%%%%%%%%%%%%%
%%%%%%%%%%%%%%% first tree
%%%%%%%%%%%%%%%%%%%%%%%%%%%%%

\put(5,55){\large $\dotsb\stackrel{\tau_{h-2}}{\longrightarrow}$}

%%%%%%%%%%%%%%% zero level
\put(125,80){\circle*{2}}
\put(80,84){$n-n_2-\dotsb-n_{h-2}-m_{h-1}$}

%%%%%%%%%%%% first level 
\put(125,80){\line(-2,-1){20}}
\put(125,80){\line(2,-1){20}}
\put(145,70){\circle*{2}}
\put(105,70){\circle*{2}}

%%%%%% labelings
\put(80,70){$n_{h-2}$}
\put(150,70){$m_{h-1}$}
\put(205,55){\large $\stackrel{\tau_{h-1}}{\longrightarrow}$}

%%%%%%%%%%%% second level
\put(105,70){\line(-2,-3){14}}
\put(105,70){\line(2,-3){14}}
\put(91,49){\circle*{2}}
\put(119,49){\circle*{2}}

\put(145,70){\line(-2,-3){14}}
\put(145,70){\line(2,-3){14}}
\put(131,49){\circle*{2}}
\put(159,49){\circle*{2}}

%%%%%%%%%%%% labelings
\put(68,49){$n_{h-3}$}
\put(102,44){$\alpha_{h-2}$}
\put(132,46){$\alpha_{h-1}$}
\put(160,46){$\alpha_h$}

%%%%%%%%%%%% third level
\put(91,49){\line(2,-3){14}}
\dottedline{4}(91,49)(77,28)
\put(77,28){\circle*{2}}
\put(105,28){\circle*{2}}

%%%%%%%%%%%% labelings
\put(64,28){$n_2$}
\put(108,26){$\alpha_{h-3}$}

%%%%%%%%%%%% last level
\put(77,28){\line(-2,-3){14}}
\put(77,28){\line(2,-3){14}}
\put(63,7){\circle*{2}}
\put(91,7){\circle*{2}}

%%%%%%%%%%%% labelings
\put(51,7){$\alpha_1$}
\put(94,7){$\alpha_2$}

%%%%%%%%%%%%%%%%%%%%%%%%%%%%
%%%%%%%%%%%%%%% second tree
%%%%%%%%%%%%%%%%%%%%%%%%%%%%%

%%%%%%%%%%%%%%% zero level
\put(300,80){\circle*{2}}
\put(290,84){$n_h$}

%%%%%%%%%%%% first level 
\put(300,80){\line(-1,-1){30}}
\put(300,80){\line(1,-1){15}}
\put(315,65){\circle*{2}}
\put(285,65){\circle*{2}}

%%%%%% labelings
\put(260,65){$n_{h-1}$}
\put(318,65){$\alpha_h$}

%%%%%%%%%%%% second level
\put(285,65){\line(1,-1){15}}
\put(270,50){\circle*{2}}
\put(300,50){\circle*{2}}

%%%%%%%%%%%% labelings
\put(243,50){$n_{h-2}$}
\put(302,50){$\alpha_{h-1}$}

%%%%%%%%%%%% third level
\put(270,50){\line(1,-1){15}}
\dottedline{4}(270,50)(255,35)
\put(255,35){\circle*{2}}
\put(285,35){\circle*{2}}

%%%%%%%%%%%% labelings
\put(288,35){$\alpha_{h-2}$}
\put(240,35){$n_2$}

%%%%%%%%%%%% last level
\put(255,35){\line(-1,-1){15}}
\put(255,35){\line(1,-1){15}}
\put(240,20){\circle*{2}}
\put(270,20){\circle*{2}}

%%%%%%%%%%%% labelings
\put(234,11){$\alpha_1$}
\put(268,11){$\alpha_2$}

\end{picture}

\noindent
It is convenient to list explicitly the coefficients of these transplantations:
\[
\begin{split}
&\text{coefficients of }\tau_1: \quad s=1,\; r=2, \;i=l=0,\; j=m_3+\dotsb+,m_{h-1},\; v=m_2+j, u=n_2;\\
&\text{coefficients of }\tau_2: \quad s=2,\; r=3, \; i=n_2,\;l=0,\; j=m_4+\dotsb+,m_{h-1},\; v=m_3+j, u=n_2+n_3;\\
&\text{coefficients of }\tau_3: \quad s=3,\; r=4, \; i=n_2+n_3,\;l=0,\; j=m_5+\dotsb+,m_{h-1},\; v=m_4+j, u=n_2+n_3+n_4;\\
&\dotsb\dotsb\dotsb\\
&\text{coefficients of }\tau_{h-1}: \quad s=h-2,\; r=h-1, \; i=n_2+\dotsb+n_{h-2},\;l=0,\; j=0,\; v=m_{h-1}, u=i+n_{h-1}.
\end{split}
\]

\noindent
Therefore from \eqref{multqRacahexpression} the resulting $(h-2)$-dimensional $q$-Racah polynomial is:

\begin{equation}\label{myGaRamultqRacah}
\begin{split}
&r_{n_2}(m_2;\alpha_2,\alpha_1,\alpha_2\dotsb\alpha_h q^{n+m_3+\dotsb+m_{h-1}+h-2},m_1+m_2|q)\\
&\times q^{-n_2m_3}r_{n_3}(m_3;\alpha_3,\alpha_1\alpha_2q^{2n_2+1},\alpha_3\dotsb\alpha_h q^{n+m_4+\dotsb+m_{h-1}-n_2+h-3},m_1+m_2+m_3-n_2|q)\\
&\times q^{-(n_2+n_3)m_4}r_{n_4}(m_4;\alpha_4,\alpha_1\alpha_2\alpha_3q^{2n_2+2n_3+2},\alpha_4\dotsb\alpha_h q^{n+m_5+\dotsb+m_{h-1}-n_2-n_3+h-4},\\
&\qquad\qquad\qquad m_1+m_2+m_3+m_4-n_2-n_3|q)\times\dotsb\dotsb\times\\
&\times q^{-(n_2+\dotsb+n_{h-2})m_{h-1}}
r_{n_{h-1}}(m_{h-1};\alpha_{h-1},\alpha_1\dotsb\alpha_{h-2}
q^{2n_2+\dotsb+2n_{h-2}+h-3},\alpha_{h-1}\alpha_h q^{n_{h-1}+n_h+1},n_{h-1}+n_h|q)\\
\equiv&\prod_{k=2}^{h-1}q^{-m_ki_k}r_{n_k}(m_k;\alpha_k,A_{k-1}q^{2i_{k-1}+k-2},A_{k-1}^{-1}A_h q^{n+j_{k+1}-i_k+h-k},n-i_k-j_k|q),
\end{split}
\end{equation}

\noindent
where $i_k=n_2+n_3+\dotsb+n_{k-1}$ and $j_k=m_{k+1}+m_{k+2}+\dotsb+m_{h-1}$. 
If we take \eqref{myGaRamultqRacah} for $n_2,n_3,\dotsc,n_h\geq 0$, $n_2+n_3+\dotsb+n_h=n$ we get a complete family of orthogonal polynomials on the set $\{(m_1,m_2,\dotsc,m_{h-1}):m_1+m_2+\dotsb+m_{h-1}=n, m_i\geq0,i=1,2,\dotsc,h-1\}$. The weight is the reciprocal of \eqref{normxi} (for $N=n$) while the square of the norm is given by the reciprocal of \eqref{normtheta} (again for $N=n$). \\

Now we consider the $q$-Racah polynomials in \cite{GaRa}. 
We use the notation of Gasper and Rahman except that we denote by $\widetilde{r}_n$ their $q$-Racah polynomial (see \eqref{GaRaqRacah}) and that we set $\widetilde{A}_k=a_1a_2\dotsb a_k$ (we recall also that $N_k=n_1+n_2+\dotsb+n_k$ and $x_{s+1}=N$). Using the symmetry property (2.2) in \cite{GaRa}, the $q$-Racah polynomials (2.6) in the same paper may be written in the form:

\begin{multline}\label{GaRamultqRacah}
R_{\bf n}({\bf x};a_1,a_2,\dotsc,a_{s+1},b,N|q)=\\
\prod_{k=1}^s\widetilde{r}_{n_k}(x_{k+1}-x_k;a_{k+1}q^{-1},b\widetilde{A}_kq^{2N_{k-1}}a_1^{-1},\widetilde{A}_k^{-1}q^{-x_{k+1}-N_{k-1}},x_{k+1}-N_{k-1}|q).
\end{multline}

\noindent
Setting 

\begin{equation}\label{transfparam}
\begin{split}
&s=h-2,\quad\qquad x_1=m_1, \qquad\qquad\qquad x_{k+1}-x_k=m_{k+1} \;\text{ for }\; k=1,2,\dotsc,h-3,\\
& b=\alpha_1,\quad\qquad a_1=\alpha_2^{-1}\alpha_3^{-1}\dotsb\alpha_h^{-1}q^{-2n-h+2}, \quad\qquad a_k=\alpha_kq \;\text{ for }\;k=2,3,\dotsc,h-1,
\end{split}
\end{equation}

\noindent
replacing $n_k$ with $n_{k+1}$ and then $k$ with $k-1$ \eqref{GaRamultqRacah} becomes:

\begin{equation}\label{GaRamultqRacah2}
\prod_{k=2}^{h-1}\widetilde{r}_{n_k}(m_k;\alpha_k,A_{k-1}q^{2i_{k-1}+k-2},A_{k-1}^{-1}A_h q^{n+j_{k+1}-i_k+h-k},n-i_k-j_k|q),
\end{equation}

\noindent
where $i_k$ and $j_k$ are as above. 
In virtue of \eqref{myGaRamultqRacah} and the conversion formula \eqref{GaRaqRacah} we get that \eqref{GaRamultqRacah2} is equal to \eqref{myGaRamultqRacah} multiplied by:

\[
(-1)^n\prod_{k=2}^{h-1}(A_kq^{2i_k+n_k+k-1},\alpha_kq,q;q)_{n_k}(A_{k-1}^{-1}A_hq^{n_k+h-k+1})^{-n_k/2}.
\]

\noindent
Similarly, using the transformations \eqref{transfparam} of the parameters, one can see (after a lot of elementary calculations) that the weight (2.16) in \cite{GaRa} is equal to our weight (the reciprocal of \eqref{normxi} for $N=n$) multiplied by:

\[
(A_{h-1}q^{h-1})^nq^{(n^2-3n)/2}\frac{(q,\alpha_hq,\alpha_{h-1}\alpha_hq^{n+1};q)_n}{(\alpha_{h-1};q)_n}.
\]

\end{document}